\documentclass[12pt]{amsart}
\usepackage{amsmath,amsfonts,amssymb,amsthm,epsfig}
\usepackage{latexsym}
\usepackage{amssymb}

\usepackage{hyperref}

\usepackage{setspace}

\usepackage{tikz}
\usetikzlibrary{shapes,arrows,positioning,decorations.markings}

\usepackage[margin=1in]{geometry}

\def\({\left(}
\def\){\right)}
\def\[{\left[}
\def\]{\right]}
\def\<{\left\langle}
\def\>{\right\rangle}

\newtheorem{theorem}{Theorem}[section]
\newtheorem{lemma}[theorem]{Lemma}
\newtheorem{definition}[theorem]{Definition}

\newtheorem{remark}[theorem]{Remark}

\DeclareMathOperator{\sign}{sign}

\newcommand{\lozl}[2]
{
	\begin{scope}[shift={#1}]
		\draw [thick,fill=#2] (0,0) -- (0,1) -- (1,0.5) -- (1,-0.5) -- cycle;
	\end{scope}
}
\newcommand{\lozr}[2]
{
	\begin{scope}[shift={#1}]
		\draw [thick,fill=#2] (0,0) -- (0,1) -- (-1,0.5) -- (-1,-0.5) -- cycle;
	\end{scope}
}

\begin{document}

\title[Turning points for periodic weights of arbitrary period]
{Turning point processes in plane partitions with periodic weights of arbitrary period}

\author[Sevak Mkrtchyan]{Sevak Mkrtchyan}
\email{sevak.mkrtchyan@rochester.edu}
\address{University of Rochester\\ Rochester, NY}

\begin{abstract}
We study random plane partitions with respect to volume measures with periodic weights of arbitrarily high period. We show that near the vertical boundary the system develops up to as many turning points as the period of the weights, and that these turning points are separated by vertical facets which can have arbitrary rational slope. In the lozenge tiling formulation of the model the facets consist of only two types of lozenges arranged in arbitrary periodic deterministic patterns. We compute the correlation functions near turning points and show that the point processes at the turning points can be described as several GUE-corners processes which are non-trivially correlated. 

The weights we study introduce a first order phase transition in the system. We compute the limiting correlation functions near this phase transition and obtain a process which is translation invariant in the vertical direction but not the horizontal.
\end{abstract}

\dedicatory{To Nicolai Reshetikhin on the occasion of his 60'th birthday.}

\maketitle


\section{Introduction}

Recall that a plane partition confined to a $c\times d$ box is a two-dimensional array $\pi=(\pi_{i,j})_{1\leq i\leq c,1\leq j\leq d}$ of non-negative integers such that $\pi_{i,j}\geq\pi_{i+1,j}$ and $\pi_{i,j}\geq\pi_{i,j+1}$ for all $i,j$. If $c$ or $d$ is infinite, then we also require that $\pi_{i,j}=0$ if $i+j$ is large enough. We will denote the set of all plane partitions confined to a $c\times d$ box by $\Pi^{c,d}$. 

A natural measure on plane partitions is the so-called ``volume'' measure, where the probability $\mathbb{P}_q(\pi)$ of a plane partition $\pi$ is proportional to $q^{|\pi|}$, where $q$ is a parameter in $(0,1)$ and $|\pi|=\sum_{i,j}\pi_{i,j}$ is the volume of $\pi$. Scaling limits of plane partitions and their generalizations under volume measures and  with various boundary conditions have been studied extensively (see e.g. \cite{NienHilBloTriangularSOS},  \cite{KO}, \cite{OR1}, \cite{OR2}, \cite{BMRT}, \cite{M} and references therein). In particular, in \cite{OR1} a broad family of measures called the Schur process was introduced and it was shown that the ``volume'' measure is a special case of the Schur process. Using the vertex operator formalism it was shown in \cite{OR1} that the Schur process is determinantal and an explicit contour-integral representation of the kernel was obtained, paving the way for studying the asymptotics of the local correlation functions for the volume measure. 

In the formulation of the model as a Schur process a plane partition is thought of as a sequence of (interlacing) Young diagrams. More precisely, given a plane partition $\pi=(\pi_{i,j})$ and an integer $-c<t<d$, let $\pi(t)$ be the $t$'th diagonal slice of $\pi$ defined by 
$$\pi(t)=
\left\{\begin{array}{ll}
(\pi_{1,k+1},\pi_{2,k+2},\pi_{3,k+3},\dots),& k\geq 0,\\
(\pi_{-k+1,1},\pi_{-k+2,2},\pi_{-k+3,3},\dots),&k\leq 0.
\end{array}\right.
$$
This decomposes $\pi$ into a sequence of interlacing Young diagrams $$\pi\Leftrightarrow \pi(-c+1)\prec\dots\prec\pi(-1)\prec\pi(0)\succ\pi(1)\succ\dots\succ\pi(d-1),$$
where the notation $\mu\succ\nu$ for Young diagrams $\mu=(\mu_1,\mu_2,\dots)$ and $\nu=(\nu_1,\nu_2,\dots)$ stands for the interlacing condition $$\mu_1\geq\nu_1\geq\mu_2\geq\nu_2\dots.$$ The volume measure can be expressed as 
$$\mathbb{P}_q(\pi)\propto\prod_{-c<t<d}q^{|\pi(t)|},$$
where $|\pi(t)|$ stands for the sum of the entries in $\pi(t)$. In this formulation it is natural to consider a generalization of the volume measure where the weight $q$ depends on the ``time'' parameter $t$ (see \cite{OR2}). Namely, given a set of positive real numbers $\bar{q}=(q_t)_{-c<i<d}$, define the measure $\mathbb{P}_{\bar{q}}$ on the space of plane partitions $\Pi^{c,d}$ by
\begin{equation}
\label{eq:qvol}
\mathbb{P}_{\bar{q}}(\pi)\propto \prod_{-c<t<d}q_t^{|\pi(t)|},
\end{equation}
i.e. give weight $q_t$ to the boxes on slice $t$ of the plane partition. The case $q_t=q,\forall t$ corresponds to the volume measure introduced earlier and we will call it the homogeneous case. The measure \eqref{eq:qvol} with inhomogeneous weights is still a Schur process. In \cite{OR2} Okounkov and Reshetikhin showed that it is a determinantal point process, gave a contour-integral representation for the correlation kernel and studied the asymptotic behavior in the case of homogeneous weights. In this paper we will analyze the scaling limit of this measure with inhomogeneous periodic weights. 

\begin{figure}[ht]
\caption{\label{fig:PlanePartition} A plane partition $\pi$ and the corresponding stack of cubes/lozenge tiling}
\includegraphics[width=9cm]{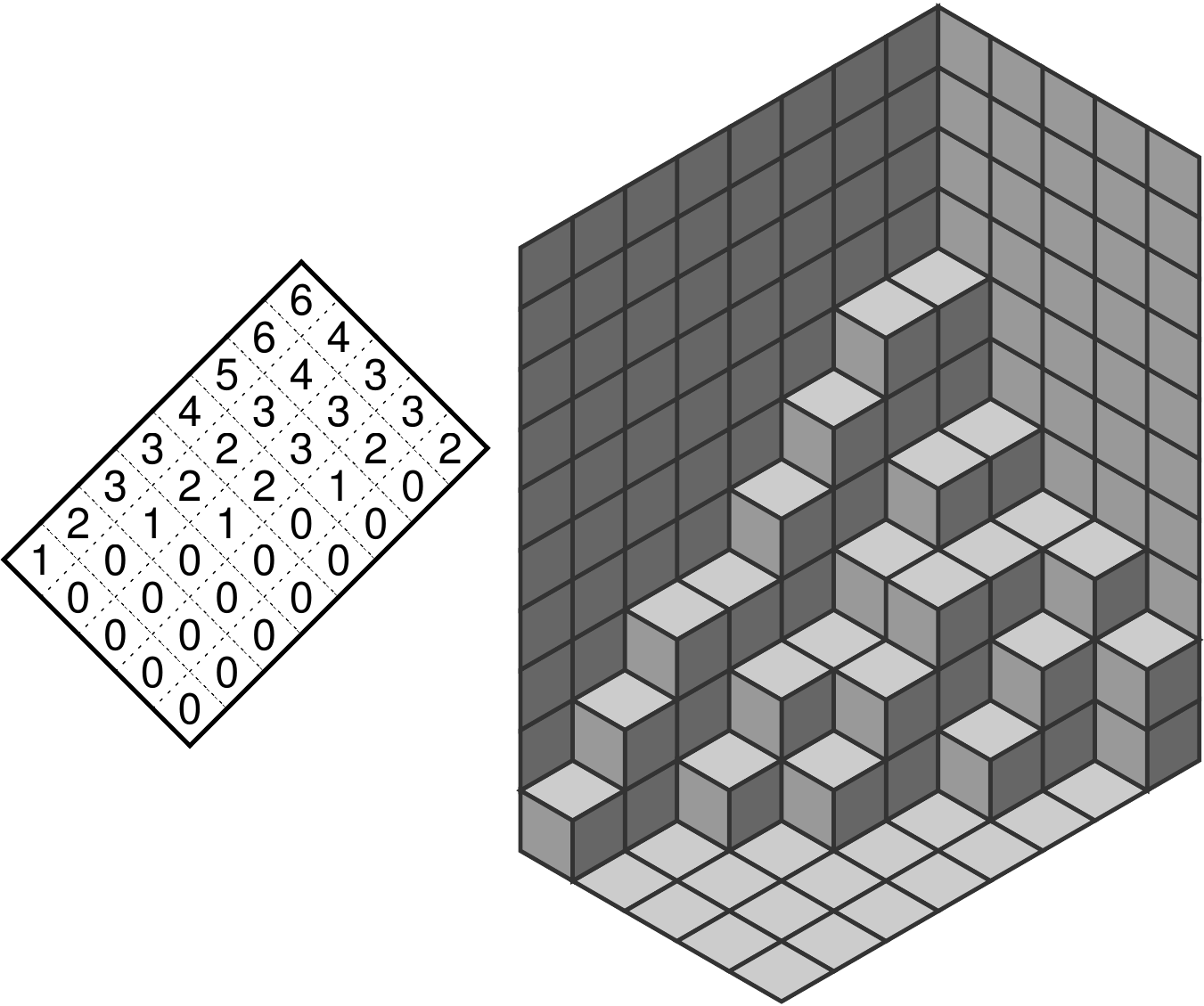}
\end{figure}

To study the scaling limit it is convenient to use a different formulation of the model. Plane partitions can be visualized through their height-functions, which can be obtained by stacking $\pi_{i,j}$ identical cubes at position $(i,j)$ for all $(i,j)$ (see Figure\ref{fig:PlanePartition}). It is evident from the figure that if we regard it as a two-dimensional object, plane partitions can be thought of as a tiling of a certain region of the plane by rhombi of three different orientations, called lozenges. 

In the appropriate scaling limit the system develops frozen and liquid regions separated by a curve called the arctic curve or the frozen boundary. In a frozen region randomness disappears and the lozenges form a deterministic pattern, whereas in a liquid region any possible local configuration of lozenges appears with positive probability. Each frozen region that forms in the scaling limit of random plane partitions distributed with respect to the homogeneous volume measure consists of only one type of lozenges, and in the 3-dimensional formulation of the model via its height function it corresponds to a facet with a normal vector pointing in one of the three lattice directions.

Along the frozen boundary there are special points called turning points. These are the points where two different frozen regions and the liquid region meet. In \cite{OR3} Okounkov and Reshetikhin conjectured that the point processes at turning points should converge to the Gaussian Unitary Ensemble corners process and proved that this is indeed the case for the homogeneous volume measure with certain boundary conditions. For results towards proving the conjecture in the case of boxed plane partitions with respect to the uniform measure see \cite{JohNordGUEminors}, \cite{GorinPanova}, \cite{Novak-GUECorners}, for Gelfand-Tsetlin patterns with the homogeneous volume measure see \cite{MP}, and for more general Gaussian orbital beta processes see \cite{Cuenca-OrbitalBetaGUE}.

\begin{figure}[ht]
\caption{\label{fig:frozenRegions} All possible frozen regions of slope $0$ that can arise when weights are $6$-periodic.}
\includegraphics[width=16cm]{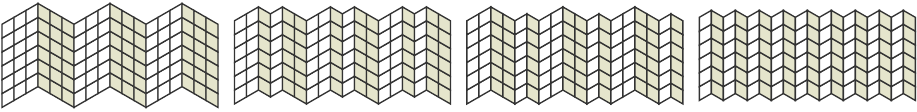}
\end{figure}

In the present paper we will mainly concentrate on the behavior of the turning points that arise in random plane partitions with respect to the non-homogeneous measure \eqref{eq:qvol} when the weights are $k$-periodic, and the new types of frozen regions that arise. We show that the point processes at turning points are no longer the GUE-corners process, but rather several copies of the GUE corners process non-trivially interlaced. Moreover, we show that the turning points are separated by new types of frozen regions. In the 3-dimensional formulation these are facets which are vertical and whose projections on the horizontal plane give a line of arbitrary rational slope. In the lozenge tiling formulation each facet consists of two types of lozenges, arranged in an arbitrary periodic deterministic pattern (see Figure \ref{fig:frozenRegions} for examples). The exact value of the slope, or the exact pattern, are determined by the specific weights taken.

\begin{figure}[ht]
\caption{\label{fig:kperfunddom} The smaller region is the fundamental domain when weights are $4$-periodic.}
\includegraphics[width=6cm]{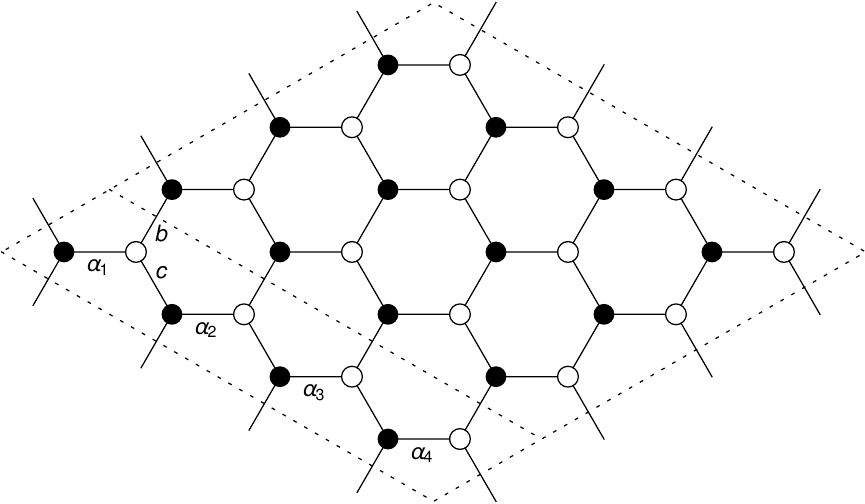}
\end{figure}

\begin{figure}[ht]
\caption{\label{fig:homogfunddom} The smaller region is the fundamental domain when weights are homogeneous.}
\includegraphics[width=6cm]{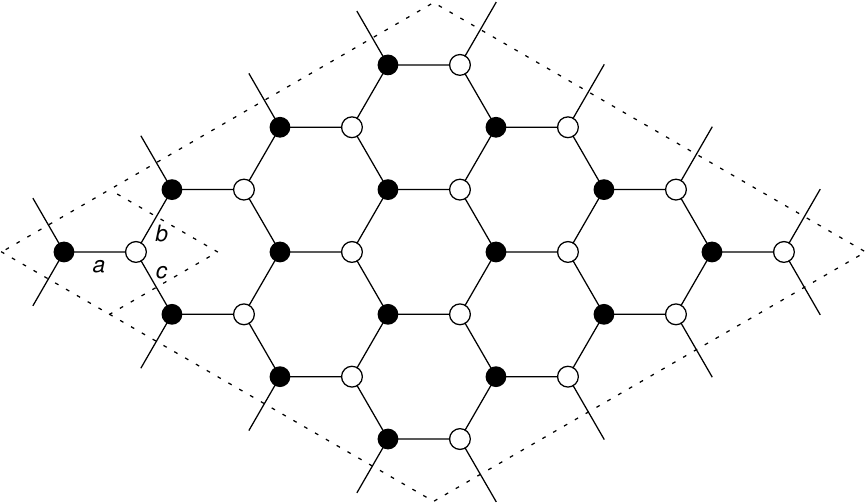}
\end{figure}

In the formulation of the model as a dimer model the weights we take essentially correspond to taking a larger fundamental domain in the language of \cite{Kenyon-DimerLectures} (see Figure \ref{fig:kperfunddom}) 
as opposed to the smaller fundamental domain from Figure \ref{fig:homogfunddom} in the case of homogeneous weights. The larger fundamental domain results in a larger Newton Polygon (see \cite{KO},\cite{KOS}), which has integer points on the boundary. As was discussed in \cite{KO}, \cite{KOS}, these integer points are responsible for the new types of frozen facets. Because the weights we take are non-homogeneous only in the horizontal direction, the Newton polygon doesn't have any integer points in its interior, which means we will not see gaseous regions in our setup.

Certain doubly periodic models have been studied in the literature. See, for example, \cite{ChhitaJohanssonTwoPer} or \cite{Berggren-Duits-Blocks}. One significant difference they have is that the systems under doubly periodic weights can develop gaseous regions as well. 

When dealing with the measure \eqref{eq:qvol} with appropriate non-homogeneous weights care should be taken to ensure that the measure is indeed well-defined, as the partition function becomes infinite. The case of $2$-periodic weights was studied in \cite{M2Per} and in that case this issue was resolved by modifying the boundary and studying certain skew plane partitions instead. This approach, however, does not work for $k>2$. Instead, what we do here is to modify the weight in one position, namely the $0$'th slice at the corner, so weights are periodic everywhere except this slice. This should have no effect on the processes at turning points, as they are macroscopically far away from that slice. However, the modification we introduce does affect the behavior of the system at the $0$'th slice. Namely, it introduces a first-order phase transition, and the local point process loses translation invariance in the horizontal direction. For a precise description of this process see Theorem \ref{thm:bulkmiddle}.

\subsection{Acknowledgments}
I would like to thank the referee for a careful reading of the paper and in particular for bringing to my attention the work of Borodin and Shlosman \cite{BorShl2009}.

The work was partially supported by the Simons Foundation Collaboration Grant No. 422190.

\section{Notation and description of main results}

\begin{figure}[ht]
\caption{\label{fig:coordinates} Coordinate axes}
\includegraphics[width=6cm]{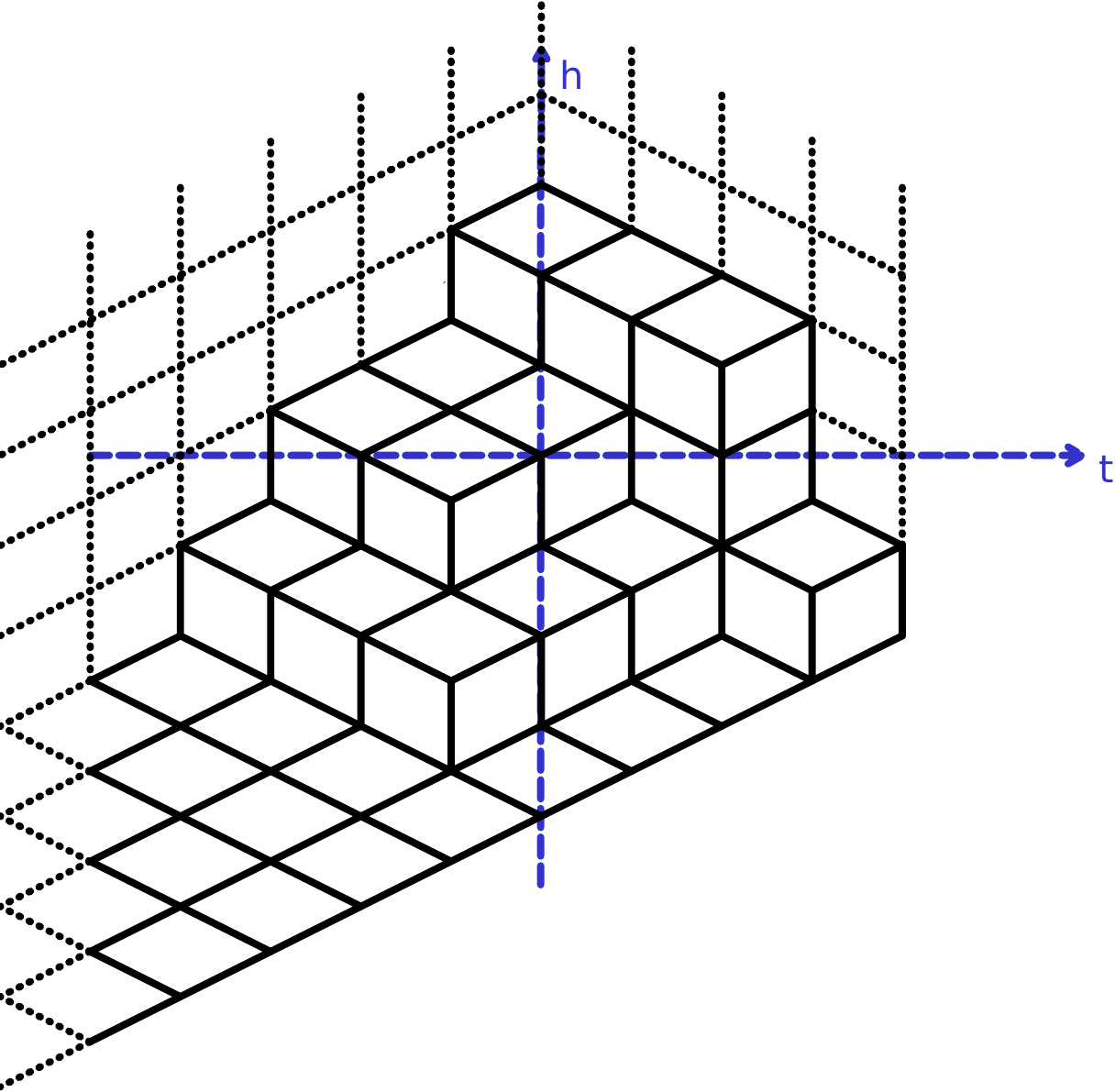}
\end{figure}

\textbf{The boundary:}
We will work with the formulation of the model in terms of lozenge tilings. Knowing the positions of lozenges of only one type completely determines the tiling. In particular, the positions of the horizontal tiles (those at the top of a stack of cubes) determine the complete tiling. Introduce planar coordinates $(t,h)$, so that the centers of horizontal tiles are on the lattice $\mathbb{Z}\times \frac 12 \mathbb{Z}$ and position the point $(0,0)$ as in Figure \ref{fig:coordinates}.

To simplify computations we will consider the case $c=\infty$, so the plane partition will consist of the slices $\pi(t)$ for $t\in(-\infty,d-1]\cap\mathbb{Z}$. 

\textbf{Weights:} We will consider a weight sequence $\bar{q}$ which is periodic of period $k\in\mathbb{N}$, i.e. $q_i=q_{i+k}$ for all $i\in\mathbb{Z}$.  The thermodynamic limit in the homogeneous case corresponds to the limit $q_t=q\rightarrow 1^-$. In the periodic case one could consider the limit when $q_i\rightarrow 1^-, \forall i$, however in that regime the behavior of the system is the same as that of the homogeneous measure with the weight given by the geometric average of the periodic weights. The more interesting limit to study is when $q_i\rightarrow\alpha_i$, where $\alpha_0,\dots,\alpha_{k-1}$ are positive numbers. When $\prod_{i=0}^{k-1}\alpha_i<1$, the expected size of the system is finite, whereas when $\prod_{i=0}^{k-1}\alpha_i>1$ the measure does not make sense as the partition function is infinite. Thus, the appropriate setup is when $\prod_{i=0}^{k-1}\alpha_i=1$. However, having any $\alpha_i>1$ still makes the partition function infinite, so the measure $P_{\bar{q}}$ is not well defined. To see this, consider cases based on the value of $\alpha_{k-1}$. 

It is convenient to set $q_i=\alpha_iq$ with $0<q<1$ and consider the limit $q\rightarrow 1^-$.

\begin{figure}[ht]
\caption{\label{fig:partfcninf1} The configuration $\pi^m$ when $\alpha_{k-1}<1$.}
\includegraphics[width=6cm]{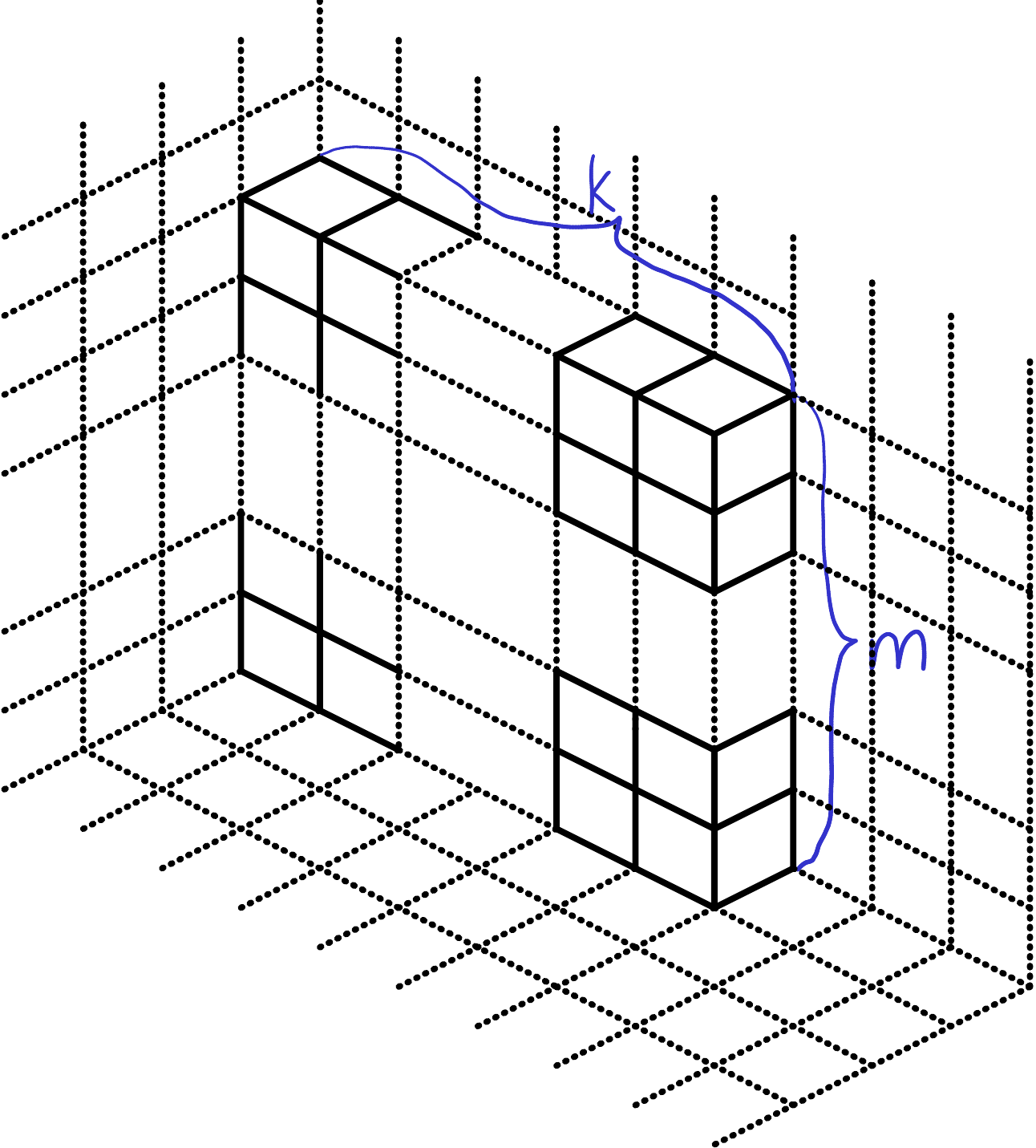}
\end{figure}

First, suppose $\alpha_{k-1}<1$. Then $\alpha_0\cdot\ldots\cdot \alpha_{k-2}>1$. Let $\pi^m$ be the plane partition given by $\pi^m_{1,j}=m$ for $1\leq j\leq k-1$ and $\pi^m_{i,j}=0$ otherwise (see Figure \ref{fig:partfcninf1}). Then $$P_{\bar{q}}(\pi^m)=(q_0\dots q_{k-2})^m\approx (\alpha_0\dots \alpha_{k-2})^m,$$ so we have $\sum_{m=1}^\infty P_{\bar{q}}(\pi^m)=\infty$, which means the measure $P_{\bar{q}}$ is not well-defined. 

\begin{figure}[ht]
\caption{\label{fig:partfcninf2} The configuration $\pi^m$ when $\alpha_{k-1}>1$.}
\includegraphics[width=6cm]{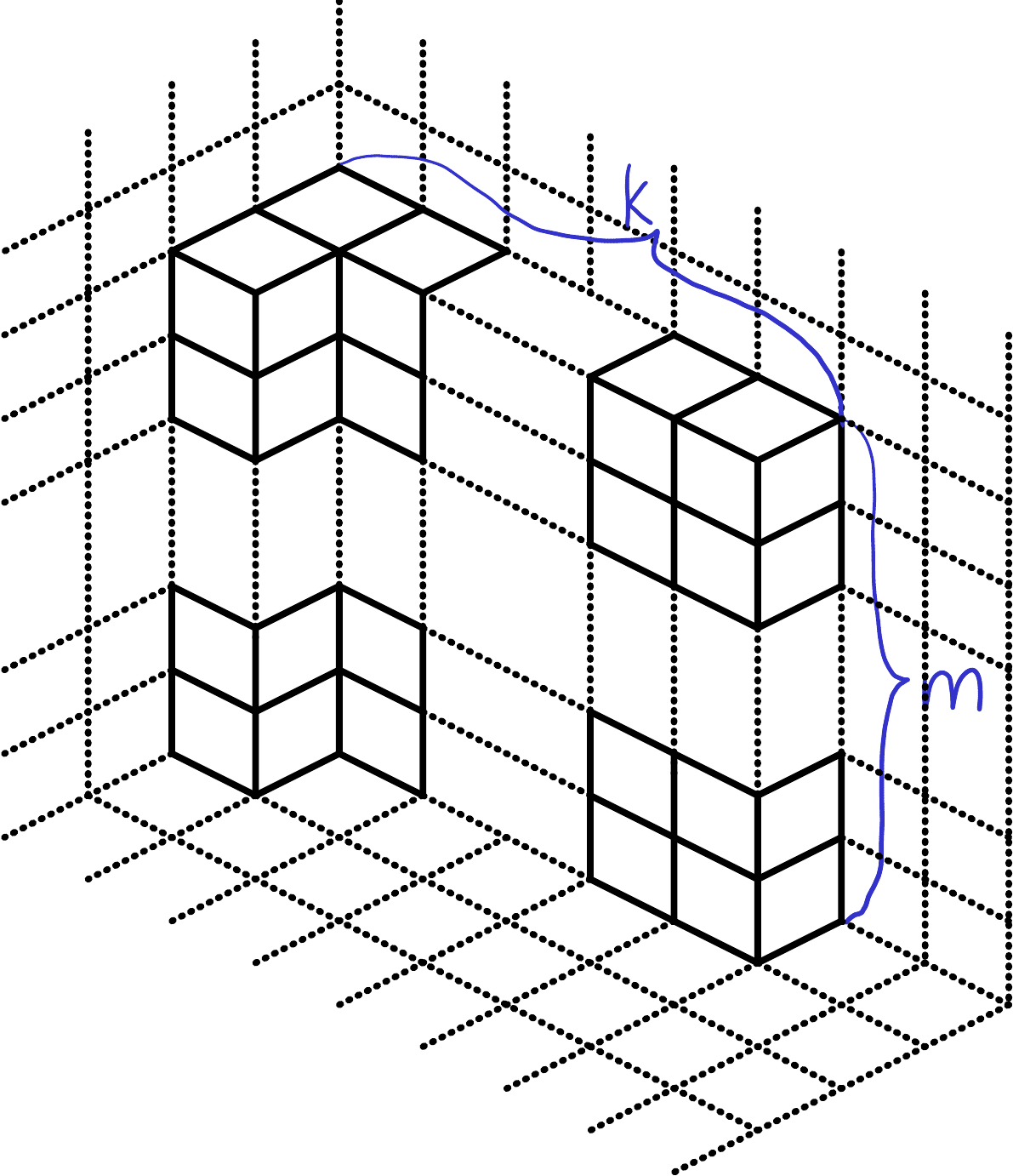}
\end{figure}

Now, suppose $\alpha_{k-1}>1$. Then $q_{-1}=q_{k-1}\approx \alpha_{k-1}>1$. Let $\pi^m$ be the plane partition given by $\pi^m_{1,j}=m$ for $1\leq j\leq k$, $\pi^m_{2,1}=m$ and $\pi^m_{i,j}=0$ otherwise (see Figure \ref{fig:partfcninf2}). Then $$P_{\bar{q}}(\pi^m)=(q_{-1}q_0\dots q_{k-1})^m\approx (\alpha_{k-1}\alpha_0\dots \alpha_{k-1})^m=\alpha_{k-1}^m,$$ so we again have $\sum_{m=1}^\infty P_{\bar{q}}(\pi^m)=\infty$. 

Lastly, if $\alpha_{k-1}=1$, then we can consider $\alpha_{k-2}$, and so on. This shows, that unless $\alpha_0=\dots=\alpha_{k-1}=1$, the measure $P_{\bar{q}}$ is not well defined if $q$ is close enough to $1$.

\begin{figure}[ht]
\caption{\label{Fig:OneCornerWeights} Weights assignments of period $4$. The number in a cell indicates the weight of the boxes above the cell.}
\includegraphics[width=7cm]{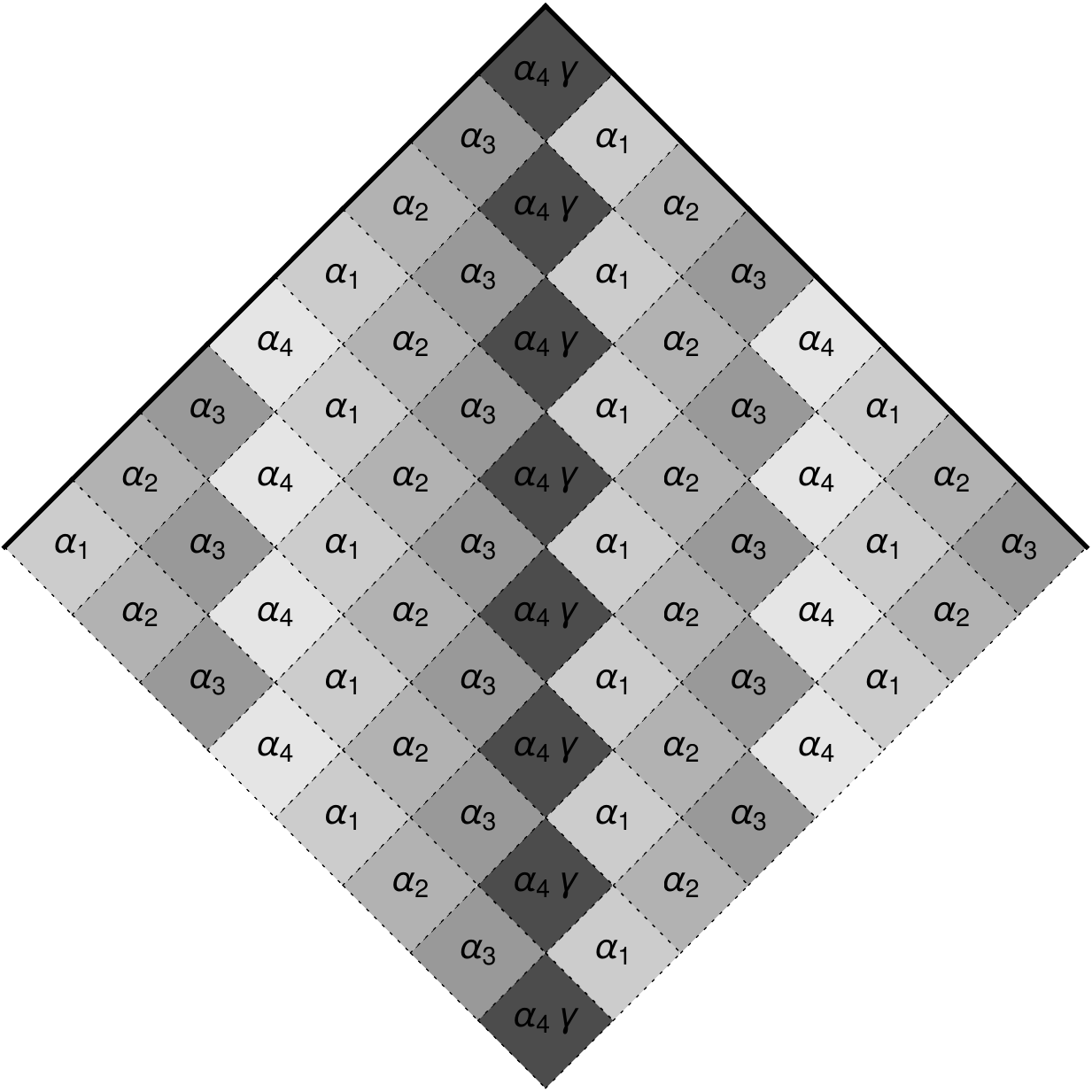}
\end{figure}

One way to fix the issue is to modify the boundary. This was the approach taken in \cite{M2Per}, where skew plane partitions with 2-periodic weights with a staircase cut-out at the corner were considered. Such an approach does not work when the period $k>2$. Instead, in this paper, we consider ordinary plane partitions, but modify the weights so they are periodic everywhere except the slice containing the corner. More precisely, let $\alpha_0,\alpha_1,\dots,\alpha_{k-1}>0$ with $$\alpha_0\cdot\dots\cdot\alpha_{k-1}=1$$ and let $$\gamma=\prod_{j:\alpha_j<1}\alpha_j.$$ Define 
\begin{equation}
\label{eq:perwts}
q_{i}:=q_{i,r}:=
\left\{\begin{array}{cc}
q\alpha_{i\bmod k},&i\neq 0\\
q\alpha_0\gamma,&i=0,
\end{array}\right.
\end{equation}
where $q=e^{-r}$. We will study random plane partitions distributed according to the measure \eqref{eq:qvol} with weights given by \eqref{eq:perwts}, in the limit when $r \rightarrow 0^+$ and $d$ grows at the rate $d=V/r$ for some constant $V>0$. To simplify formulas we will assume that $d\bmod k$ is independent of $r$. For an illustration of the weight assignments see Figure \ref{Fig:OneCornerWeights}.

Note, that the modification of the weight $q_0$ by $\gamma$ will not have any effect on the processes observed at turning points since the turning points will be macroscopically away from the $0$'th slice. Indeed, when $k=2$, the process observed here matches with that obtained in \cite{M2Per}, where the weights are completely 2-periodic. The addition of $\gamma$ to the slice at zero does, however, has an effect on the local processes near that slice. Whereas the processes appearing in the bulk are translation invariant in two directions, the addition of $\gamma$ causes a first order phase transition with a local point process only translation invariant in the vertical direction near the special slice. 

All fully $\mathbb{Z}\times\mathbb{Z}$ translation invariant ergodic Gibbs measures were classified in \cite{SheffieldThesis}. The process we obtain away from the phase transition is translation invariant under translations by elements of $k\mathbb{Z}\times \mathbb{Z}$ only. As a result it is not given by the incomplete-beta kernel, but rather a deformation of it, which is a special case of processes studied previously in \cite{Bor} and \cite{BorShl2009}.

\textbf{Main results:} We will study the scaling limit of our model via the local correlation functions. Given a subset $U=\{(t_1,h_1),\ldots,(t_n,h_n)\}\subset \mathbb{Z}\times \frac 12 \mathbb{Z}$, define the corresponding local correlation function $\rho_{\bar{q}}(U)$ as the probability for a random tiling taken from the above probability space to have horizontal tiles centered at all positions $(t_i,h_i)_{i=1}^n$. 

The main theorems of the paper describe the point processes that appear in the thermodynamic limit of the system in terms of their correlation functions. Theorem \ref{thm:bulkNearBoundary} gives the limiting point process in the bulk near the edge, Theorem \ref{thm:bulkmiddle} gives the non-translation-invariant point process in the bulk along the slice at the corner where we observe a first-order phase transition, and Theorem \ref{thm:turningPoints} gives the point processes near turning points. Section \ref{sec:GUE-corners} explains the connection with the GUE-corners process and Section \ref{sec:FrozenRegions} describes the new frozen regions that arise.

\section{The correlation kernel and its leading asymptotics}

Okounkov and Reshetikhin \cite{OR2} showed that for arbitrary $\bar{q}$, the correlation functions of the Schur process are determinantal:
\begin{theorem}[Theorem 2, part 3 \cite{OR2}] 
\label{thm:fin-corr2}
The correlation functions $\rho_{\bar{q}}$ are determinants
\begin{equation*}
\rho_{\bar{q}}(U)
=\det(K_{\bar{q}}((t_i,h_i),(t_j,h_j)))_{1\leq i,j\leq n},
\end{equation*}
where the correlation kernel $K_{\bar{q}}$ is given by the double integral
\begin{multline}
\label{eq:main-corr2}
K_{\bar{q}}((t_1,h_1),(t_2,h_2))
=\\ \frac{1}{(2\pi \mathfrak{i})^2}
\int_{z\in C_z}\int_{w\in C_w}
\frac{\Phi_{\bar{q}}(z,t_1)}{\Phi_{\bar{q}}(w,t_2)}
\frac{1}{z-w}z^{-h_1-\frac 12 |t_1|+\frac 12}w^{h_2+\frac 12 |t_2|+\frac12}\frac{dz\ dw}{zw},
\end{multline}
where 
\begin{align}
\nonumber\Phi_{\bar{q}}(z,t)&=\frac{\Phi^-_{\bar{q}}(z,t)}{\Phi^+_{\bar{q}}(z,t)},\\
\label{eq:Phis}\Phi^+_{\bar{q}}(z,t)&=\prod_{\max(0,t)<m<d, m\in \mathbb{Z}+\frac 12}(1-zx_m^+),\\
\nonumber\Phi^-_{\bar{q}}(z,t)&=\prod_{m<\min(0,t), m\in \mathbb{Z}+\frac 12}(1-z^{-1}x_m^-),
\end{align}
the parameters $x^\pm_m$ satisfy the conditions
\begin{align}
\nonumber\frac{x^+_{m+1}}{x^+_m}&=q_{m+\frac 12},\ 0<m<d-1,
\\x^-_{-\frac 12}x^+_{+\frac 12}&=q_{0},
\label{eq:xs-qs}
\\\nonumber\frac{x^-_m}{x^-_{m+1}}&=q_{m+\frac 12},\ m<-1,
\end{align}
and $C_z$ (resp. $C_w$) is a simple positively oriented contour around 0 such that its interior contains none of the poles of $\Phi_{\bar{q}}(\cdot,t_1)$ (resp. all of $\Phi_{\bar{q}}(\cdot,t_2)^{-1}$). Moreover, if $t_1< t_2$, then $C_z$ is contained in the interior of $C_w$, and otherwise, $C_w$ is contained in the interior of $C_z$. 
\end{theorem}


The conditions \eqref{eq:xs-qs} can be rewritten in the form
\begin{align}
\label{eq:xs-prodqs}
x^-_m&=a^{-1}q_{m+\frac 12}q_{m+\frac 32}\cdot\ldots\cdot q_{d-1},
\\\nonumber x^+_m&=aq_{m+\frac 12}^{-1}q_{m+\frac 32}^{-1}\cdot\ldots\cdot q_{d-1}^{-1},
\end{align}
where $a>0$ is an arbitrary parameter. We will choose $a=e^{-V}=e^{-rd}$.

\subsubsection{The asymptotically leading term in the integral for the correlation kernel}
It can be shown that as $r\rightarrow 0$, the volume of a typical plane partition grows at the rate of $1/r^3$ so we scale the system by $r$ in all directions. Let $(\tau,\chi)$ be the scaled coordinates. First, given a point $(\tau,\chi)$, we will study the asymptotics of the correlation kernel \eqref{eq:main-corr2} when $r\rightarrow 0$, $rt_i\rightarrow\tau$, $rh_i\rightarrow\chi$ and $\Delta t=t_1-t_2,\Delta h=h_1-h_2$ are fixed.

\begin{definition}
\label{def:betas}
For $i\geq 1$ let $\beta_i=\alpha_{d-1}\alpha_{d-2}\cdot\ldots\cdot\alpha_{d-i}$. Let $l$ be the number of distinct $\beta$'s, and let $\tilde{\beta}_1<\dots<\tilde{\beta}_l$ be the ordered list of distinct $\beta$'s. For $i=1,\dots,l$ let $m_i$ be the number of times $\tilde{\beta}_i$ appears among $\beta_1,\dots,\beta_k$. We have $m_1+\dots+m_l=k$.
\end{definition}
Note, that since $d\bmod k$ is constant, the values $\beta_1,\beta_2,\dots$ do not depend on $r$. Moreover, since $a_0\cdot\ldots\cdot a_{k-1}=1$, the sequence $\beta_1,\beta_2,\dots$ is also periodic of period $k$. We might, however, have some of the numbers $\beta_1,\dots,\beta_k$ be equal, so $l$ can be less than $k$ and $m_i>1$.

We have
$$q_{m+\frac 12}\cdot\ldots\cdot q_{d-1}=q^{d-m-\frac 12}\beta_{(d-m-\frac 12) \bmod k}\gamma^{c_m},$$
where $c_m=0$ if $m>0$ and $c_m=1$ if $m<0$. Thus, we have
\begin{align*}
\lim_{r\rightarrow 0}r\ln\Phi_{\bar{q}}(z,t)
=&\lim_{r\rightarrow 0}\left(-\sum_{i=1}^k
\sum_{\substack{\max(0,t)<m<d\\m\in\mathbb{Z}+\frac 12\\\left(d-m-\frac 12\right) \bmod k = i}}r\ln(1-zae^{r(d-m-\frac 12)}\beta_i^{-1})
\right.\\&\qquad\qquad+\left.
\sum_{i=1}^k\sum_{\substack{m<\min(0,t)\\m\in\mathbb{Z}+\frac 12\\\left(d-m-\frac 12\right) \bmod k = i}}r\ln(1-z^{-1}a^{-1}e^{-r(d-m-\frac 12)}\beta_i\gamma)\right)
\\=&\frac 1k\sum_{i=1}^k\left(\int_{-\infty}^{\min(\tau,0)}\ln(1-z^{-1}\beta_i\gamma e^{M})dM
-\int_{\max(\tau,0)}^V\ln(1-z\beta_i^{-1}e^{-M})dM\right),
\end{align*}
where in the last step we replaced a sum of $2k$ Riemann sums by the sum of the corresponding integrals.

It follows that the integrand in formula \eqref{eq:main-corr2} for the correlation kernel has the leading asymptotics given by
\begin{equation}
\label{eq:kernelViaS}
K_{\bar{q}}((t_1,h_1),(t_2,h_2))
=\\ \frac{1}{(2\pi \mathfrak{i})^2}
\int_{z\in C_z}\int_{w\in C_w}e^{\frac 1r (S_{\tau,\chi}(z)-S_{\tau,\chi}(w))+O(1)}\frac{dzdw}{z-w},
\end{equation}
where
\begin{multline}
\label{eq:S}
S_{\tau,\chi}(z)=\frac 1k\sum_{i=1}^k\left(\int_{-\infty}^{\min(\tau,0)}\ln(1-z^{-1}\beta_i\gamma e^{M})dM
-\int_{\max(\tau,0)}^V\ln(1-z\beta_i^{-1}e^{-M})dM\right)\\-\left(\chi+\frac 12|\tau|\right)\ln z.
\end{multline}

\section{The point processes in the bulk near the edge and at the turning points}
\label{sec:bulkandturningpoints}
To compute the asymptotics of the correlation kernel, we use the steepest descent method, so we need to study the critical points of $S_{\tau,\chi}$. We have
\begin{equation*}
z\frac{dS_{\tau,\chi}}{dz}=-\frac 1k\sum_{i=1}^k\ln\left(\frac{z-\beta_i\gamma e^{\min(\tau,0)}}{z}\right)+\frac1k\sum_{i=1}^k\ln\left(\frac{z-\beta_ie^{V}}{z-\beta_ie^{\max(\tau,0)}}\right)-V-\chi+\frac\tau2
\end{equation*}
and
\begin{equation*}
\frac{d}{dz}\left(z\frac{d S_{\tau,\chi}(z)}{dz}\right)
=-\frac1k\sum_{i=1}^k\left(\frac1{z-e^{\min(\tau,0)}\beta_i\gamma}-\frac1{z}\right)
+\frac1k\sum_{i=1}^k\left(\frac1{z-\beta_ie^V}-\frac1{z-\beta_ie^{\max(\tau,0)}}\right).
\end{equation*}

We are interested in the behavior of the system near the turning points which appear when $\tau=V$, so we set $\tau=V-\varepsilon$ with $\varepsilon>0$ and small. We need $\varepsilon$ to be smaller than a constant that only depends on $d$ and the weights $\alpha_0,\dots,\alpha_{k-1}$. In particular we have $\tau\geq 0$, so 
\begin{equation}
\label{eq:Sp}
z\frac{dS_{\tau,\chi}}{dz}(z)=-\frac 1k\sum_{i=1}^lm_i\ln\left(\frac{z-\tilde\beta_i\gamma }{z}\right)+\frac1k\sum_{i=1}^lm_i\ln\left(\frac{z-\tilde{\beta}_ie^{V}}{z-\tilde\beta_ie^{\tau}}\right)-V-\chi+\frac\tau2
\end{equation}
and
\begin{equation}
\label{eq:Spp}
\frac{d}{dz}\left(z\frac{d S_{\tau,\chi}(z)}{dz}\right)
=-\frac1k\sum_{i=1}^lm_i\left(\frac1{z-\tilde\beta_i\gamma}-\frac1{z}\right)
+\frac1k\sum_{i=1}^lm_i\left(\frac1{z-\tilde\beta_ie^V}-\frac1{z-\tilde\beta_ie^{\tau}}\right),
\end{equation}
where we have grouped terms with equal $\beta$'s together. To simplify notation, we will denote $Sp(z):=z\frac{dS_{\tau,\chi}}{dz}(z)$. 

\subsubsection{The possible locations of double real critical points}
First, we look for double real critical points of $S_{\tau,\chi}$. If $z$ is such a point, then from \eqref{eq:S} or \eqref{eq:Sp} it is easy to see that $z\notin[0,\tilde\beta_i\gamma]$ and $z\notin[\tilde\beta_ie^{\tau},\tilde\beta_ie^V]$ for any $i$. 

\begin{figure}[ht]
\caption{\label{fig:excludedintervals} The integration contours and the intervals where there cannot be real critical points.}
\includegraphics[width=14cm]{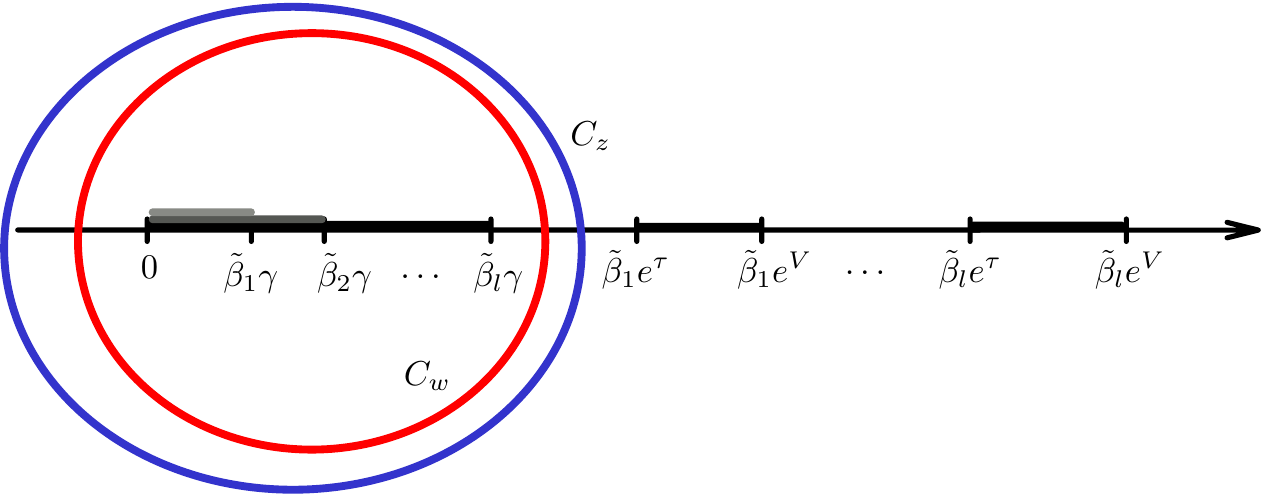}
\end{figure}

Since the $\tilde\beta$'s are distinct, by picking the bound on $\varepsilon$ to be small enough, we will have that the intervals $[\tilde\beta_ie^{\tau},\tilde\beta_ie^V]$ are disjoint. Also we have 
\begin{equation}
\label{eq:order}
\tilde\beta_i\gamma<\tilde\beta_l\gamma<\tilde\beta_1e^\tau< \tilde\beta_ie^\tau.
\end{equation}
Here the only non-trivial inequality is the middle one. It holds, since $e^\tau\geq 1$, and $\tilde\beta_1/\tilde\beta_l$ is either the product of consecutive $\alpha$'s, or the inverse of the product of consecutive $\alpha$'s. Since the product of all the $\alpha$'s is one, we have that $\tilde\beta_1/\tilde\beta_l$ is always the product of consecutive $\alpha$'s. However, $\gamma$ is the product of all $\alpha$'s which are less than $1$, so we have $\gamma<\tilde\beta_1/\tilde\beta_l$. It follows, that the intervals $[0,\tilde\beta_i\gamma]$ and $[\tilde\beta_ie^{\tau},\tilde\beta_ie^V]$ are situated on the real line as in Figure \ref{fig:excludedintervals}.

Note also, that the $C_z$ contour from \eqref{eq:main-corr2} contains no zeros of $\Phi^+_{\bar q}(z,t_1)$, i.e. none of $e^{r(m+\frac 12)}\tilde\beta_i$ for $t_1<m<d$, and the $C_w$ contour from \eqref{eq:main-corr2} contains all zeros of $\Phi^-_{\bar q}(z,t_1)$, i.e. all of $e^{r(m+\frac 12)}\tilde\beta_i\gamma$ for $m<0$. This means that the contours $C_z$ and $C_w$ cross the positive part of the real line between $\tilde\beta_l\gamma$ and $\tilde\beta_1e^\tau$ as in Figure \ref{fig:excludedintervals}.

\subsubsection{Computing the double real critical points}
It follows from the computations in the previous section that if $z$ is a double real critical point of $S_{\tau,\chi}$, then the summands in the first sum in \eqref{eq:Spp} are all of the same sign, whereas the summands in the second sum in \eqref{eq:Spp} are all small (they converge to $0$ as $\varepsilon$ goes to $0$), unless $z$ is close to $\tilde\beta_ie^\tau$ or $\tilde\beta_ie^V$, so if $z$ is a double real critical point, $z$ must be close to $\tilde\beta_ie^\tau$ or $\tilde\beta_ie^V$ for some $i$. Let $z=\tilde\beta_je^V+\delta$, where $\delta\rightarrow 0$ as $\varepsilon\rightarrow0$. All the summands in the second sum in \eqref{eq:Spp} will be of size $o(1)$ except the summand with $i=j$. The latter one should be of order $1$. Thus we have
$$\frac{\tilde\beta_j(e^V-e^\tau)}{(z-\tilde\beta_je^V)(z-\tilde\beta_je^\tau)}=\theta(1),$$
where $\theta(1)$ means of constant order as $\varepsilon\rightarrow 0$. Plugging in the values of $\tau$ and $z$, we have 
$$\frac{\tilde\beta_je^V(\varepsilon+O(\varepsilon))}{\delta(\delta+\tilde\beta_je^V\varepsilon+O(\varepsilon^2))}=\theta(1),$$
which implies that $\delta=\theta(\sqrt\varepsilon)$. From \eqref{eq:Spp} we must have
\begin{equation*}
0
=-\frac1k\sum_{i=1}^l\frac{m_i\tilde\beta_i\gamma}{\tilde\beta_je^V(\tilde\beta_je^V-\tilde\beta_i\gamma)}
+O(\delta)+O(\varepsilon)
+\frac1k\frac{m_j\tilde\beta_je^V\varepsilon}{\delta^2+O(\varepsilon^{3/2})}.
\end{equation*}
Solving for $\delta$, we get
\begin{equation}
\delta=\pm\frac{\sqrt{m_j}\tilde\beta_je^V}{\sqrt{\sum_{i=1}^l\frac{m_i\tilde\beta_i\gamma}{\tilde\beta_je^V-\tilde\beta_i\gamma}}}\sqrt{\varepsilon}+O(\varepsilon).
\end{equation}
Note, that \eqref{eq:order} implies the summands in the denominator are all positive, so $\delta$ is real. Setting $$f_j:=\frac{\sqrt{m_j}\tilde\beta_je^V}{\sqrt{\sum_{i=1}^l\frac{m_i\tilde\beta_i\gamma}{\tilde\beta_je^V-\tilde\beta_i\gamma}}},$$
let $z_{j,\pm}=\tilde\beta_je^V\pm f_j\sqrt{\varepsilon}+g_\varepsilon \varepsilon+o(\varepsilon)$. 
By setting $Sp(z_{j,\pm})=0$ and solving for $\chi$ from \eqref{eq:Sp} we obtain
\begin{equation}
\label{eq:chiCr}
\chi_{j,\pm}=-\frac V2-\frac1k\sum_{i=1}^lm_i\ln\left(\frac{e^V\tilde\beta_j-\tilde\beta_i\gamma}{e^V\tilde\beta_j}\right)\mp\sum_{i=1}^l f_j\frac{\tilde\beta_i\gamma}{e^V\tilde\beta_j(e^V\tilde\beta_j-\tilde\beta_i\gamma)}\sqrt{\varepsilon}\mp\frac1km_j\frac{e^V\tilde\beta_j}{f_j}\sqrt{\varepsilon}+O(\varepsilon).
\end{equation}
Note, that to evaluate $\chi_{j,\pm}$ up to order $\sqrt{\varepsilon}$ we do need an expansion of $z_{j,\pm}$ up to order $\varepsilon$, even though the coefficient $g_\varepsilon$ of $\varepsilon$ in $z_{j,\pm}$ cancels out in the computation of $\chi_{j,\pm}$. 

Also note, that we have
$$
z_{l,+}>z_{l,-}>z_{l-1,+}>z_{l-1,-}>\dots>z_{2,+}>z_{2,-1}>z_{1,+}>z_{1,-}
$$
and
$$
\chi_{l,+}<\chi_{l,-}<\chi_{l-1,+}<\chi_{l-1,-}<\dots<\chi_{2,+}<\chi_{2,-1}<\chi_{1,+}<\chi_{1,-}.
$$
Now, fix $\chi$ and consider the function $Sp(z)$. We have $Sp(z)<0$ if $z<0$. 

\begin{figure}[ht]
\caption{\label{fig:SpGraph} A plot of the function $Sp(z)$.}
\includegraphics[width=12cm]{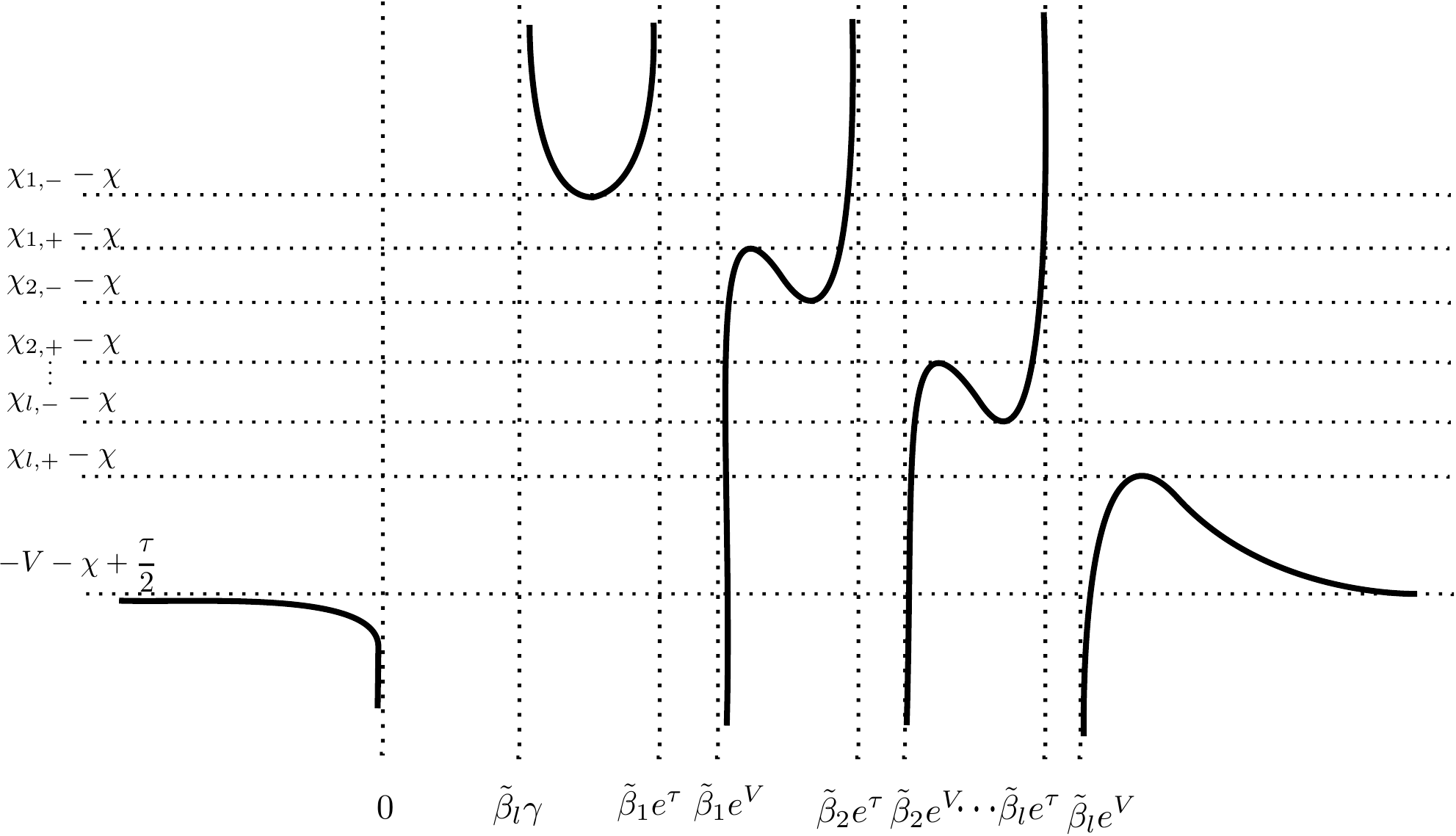}
\end{figure}

It is easy to see that 
\begin{itemize}
\item when $z$ ranges from $-\infty$ to $0$, the function $Sp(z)$ decreases from $-V-\chi+\frac \tau2$ to $-\infty$,
\item when $z$ ranges from $\tilde\beta_l\gamma$ to $\tilde\beta_1e^\tau$, the function $Sp(z)$ decreases from $-\infty$ to $\chi_{1,-}-\chi$ at $z_{1,-}$ and increases to $\infty$,
\item when $z$ ranges from $\tilde\beta_ie^V$ to $\tilde\beta_{i+1}e^\tau$, the function $Sp(z)$ increases from $-\infty$ to $\chi_{i,+}-\chi$ at $z_{i,+}$, decreases to $\chi_{i+1,-}-\chi$ at $z_{i+1,-}$ and increases to $\infty$,
\item when $z$ ranges from $\tilde\beta_le^V$ to $\infty$, the function $Sp(z)$ increases from $-\infty$ to $\chi_{l,+}-\chi$ at $z_{l,+}$ and decreases to $-V-\chi+\frac\tau2$.
\end{itemize}
These properties of $Sp(z)$ are visualized in Figure \ref{fig:SpGraph}.

\subsubsection{The number of real and complex critical points}
\label{sec:numCritPts}
First, we count the critical points of $S_{\tau,\chi}$ when $\chi$ is large. We use the following lemma proven in \cite{M}:
\begin{lemma}[Lemma 2.3 in \cite{M}]
\label{lem:geom}
Let $\{x_i\}_{i=1}^{m_1}$, $\{y_i\}_{i=1}^{m_2}$, $w$ and $\mathfrak{w}\neq 0$ be real numbers such that
$$x_{m_1}<\ldots<x_2<x_1<w<y_1<y_2<\ldots < y_{m_2}.$$
Let $\zeta$ be a complex number with $\Im{\zeta}\geq 0$. Define
\begin{equation*}
\begin{array}{c}
\mu_i=\mathrm{angle}(\zeta-x_{i+1},\zeta-x_i),\ \forall i=1,2,\ldots,m_1-1,\\
\nu_i=\mathrm{angle}(\zeta-y_i,\zeta-y_{i+1}),\ \forall i=1,2,\ldots,m_2-1,
\end{array}
\end{equation*}
and
$$
\alpha=\mathrm{angle}(\zeta-x_1,\zeta-w).
$$
\begin{figure}[ht]
\caption{\label{Fig:geom1} Setup of Lemma \ref{lem:geom}.}
\includegraphics[width=12cm]{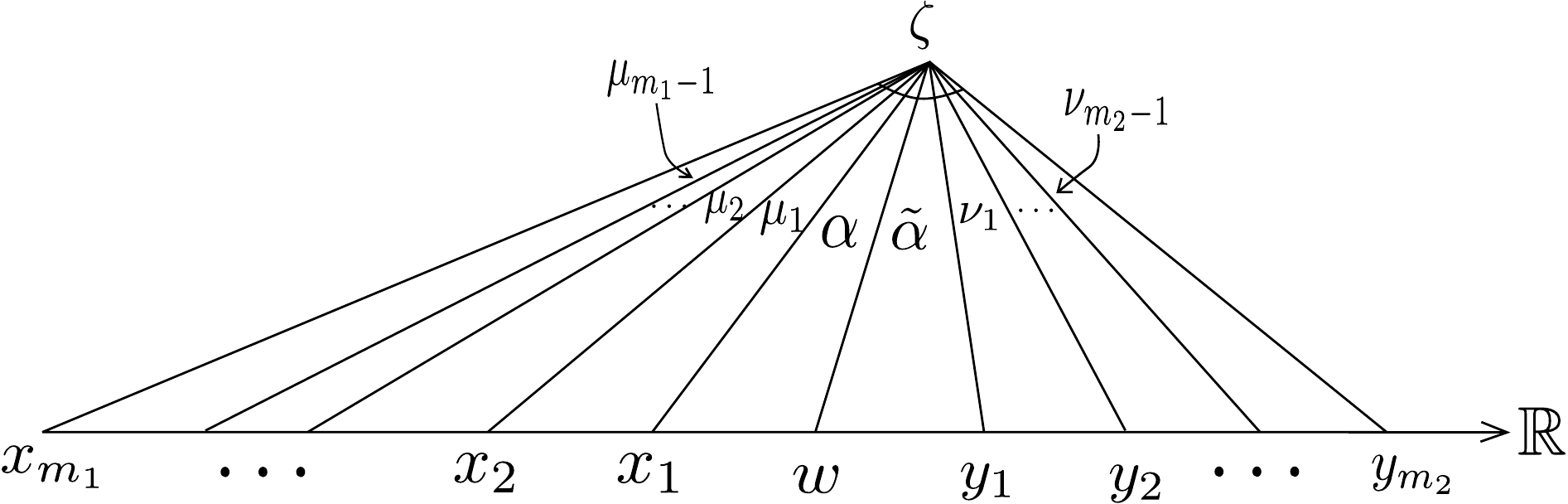}
\end{figure}
(See Figure \ref{Fig:geom1} for an illustration of the setup on the complex plane.) Given any real numbers $\{\mathfrak{x}_i\}_{i=1}^{m_1-1}$, $\{\mathfrak{y}_i\}_{i=1}^{m_2-1}$, there exists $\varepsilon>0$ such that if $|\zeta-w|<\varepsilon$, then
\begin{equation}
\label{eq:angles}
\sum_{i=1}^{m_1-1} \mathfrak{x}_i\mu_i + \mathfrak{w}\alpha +\sum_{i=1}^{m_2-1} \mathfrak{y}_i\nu_i=0
\end{equation}
if and only if $\zeta\in\mathbb{R}$ and $\zeta>w$.

The same holds if $\alpha$ is replaced by $\tilde{\alpha}:=\mathrm{angle}(\zeta-w,\zeta-y_1)$ and $\zeta>w$ is replaced by $\zeta<w$.
\end{lemma}
It follows from \eqref{eq:Sp} that if $\chi$ is large, then $z$ is close to either $\tilde\beta_l\gamma$ or $\tilde\beta_ie^\tau$ for some $i=1,\dots,l$. Setting the imaginary part of $Sp(z)$ to zero, it follows from Lemma \ref{lem:geom} that all the critical points of $S_{\tau,\chi}$ are real if $\chi$ is large. For example, if $z$ is close to $\tilde\beta_l\gamma$, then setting $\zeta=z$, $m_1=l$, $x_{m_1}=0$, $x_i=\tilde\beta_{l-i}\gamma$, $\mathfrak{x}_i=m_{l-i}+\dots+m_{l}$ for $i=1,\dots,l-1$, $w=\tilde\beta_l\gamma$, $\mathfrak{w}=m_1$,  $m_2=2l$, $y_{2j-1}=\tilde\beta_je^\tau$, $y_{2j}=\tilde\beta_je^V$, $\mathfrak{y}_{2j-1}=m_j$ and $\mathfrak{y}_{2j}=0$ for $j=1,\dots,l$ in Lemma \ref{lem:geom} gives us that $z$ must be real. 

Thus, when $\chi$ is large enough, $S_{\tau,\chi}$ has exactly $l+1$ real critical points and no non-real complex critical points. As $\chi$ varies, the number of real or complex, non-real critical points of $S_{\tau,\chi}$ can only change when it has a double real critical point. It follows from Figure \ref{fig:SpGraph} that if $\chi_{i,+}<\chi<\chi_{i,-}$ for some $i=1,\dots,l$, then $S_{\tau,\chi}$ has exactly two non-real complex critical points. Otherwise, it has only real critical points. Let $z_{cr,i}$ be the critical point with positive imaginary part when $\chi_{i,+}<\chi<\chi_{i,-}$. We moreover have that as $\chi$ increases from $\chi_{i,+}$ to $\chi_{i,-}$, the real part of $z_{cr}$ decreases from $\tilde\beta_ie^V+f_i\sqrt{\varepsilon}$ to $\tilde\beta_ie^V-f_i\sqrt{\varepsilon}$.

\subsubsection{Deformation of contours and the limiting correlation kernel in the bulk} 
\begin{figure}[ht]
\caption{\label{fig:contoursBulk} The shaded region is the region $\Re(S_{\tau,\chi}(z))>\Re(S_{\tau,\chi}(z_{cr}))$. The solid curves are the original contours whereas the dotted ones are the new contours.}
\includegraphics[width=16cm]{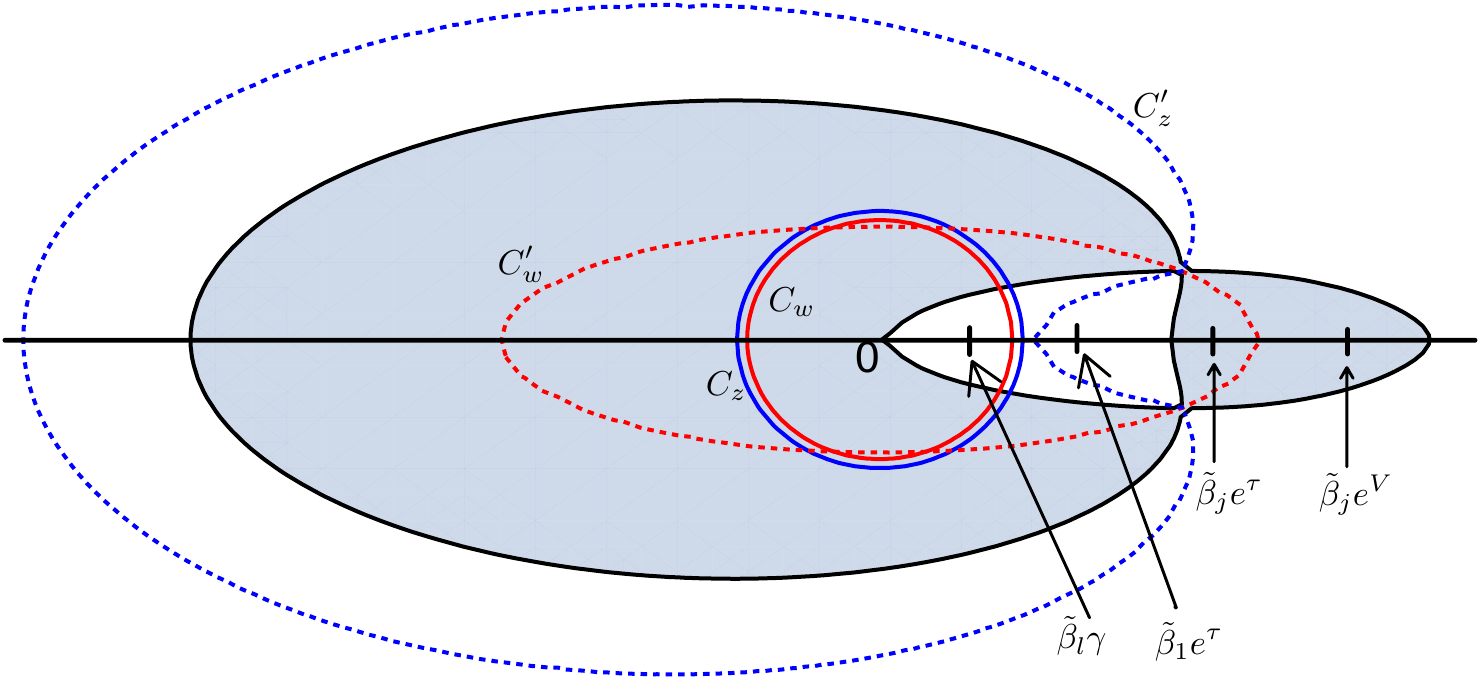}
\end{figure}

We carry out the steepest descent analysis by deforming the contours of integration in \eqref{eq:main-corr2} appropriately. Let $\chi_{i,+}<\chi<\chi_{i,-}$ and $\tau=V-\varepsilon$ with $\varepsilon$ small. Then, as we saw, the function $S_{\tau,\chi}$ has exactly one pair of complex, conjugate critical points. Let $z_{cr}$ be the one with positive imaginary part. It follows that the curve $\Re S_{\tau,\chi}(z)=\Re S_{\tau,\chi}(z_{cr})$ looks as in Figure \ref{fig:contoursBulk}. We deform the contours $C_z$ and $C_w$ to contours $C'_z$ and $C'_w$ such that they pass through $z_{cr},\bar z_{cr}$ and $\Re(S_{\tau,\chi}(z))\leq \Re(S_{\tau,\chi}(z_{cr}))\leq \Re(S_{\tau,\chi}(w))$ for all $z\in C'_z$ and $w\in C'_w$, with equality only at the critical points.

This deformation of contours results in the $z$ and $w$ contours crossing, which, together with the presence of the $\frac{1}{z-w}$ term in \eqref{eq:main-corr2}, means that as a result of the deformation we pick up residues where the curves cross, i.e. along an arc from $\bar z_{cr}$ to $z_{cr}$ which crosses the real axes to the right of $0$ if $t_1\geq t_2$ and to the left of $0$ if $t_1<t_2$. We have
\begin{multline*}
K_{\bar{q}}((t_1,h_1),(t_2,h_2))
=
\int_{z\in C'_z}\int_{w\in C'_w}e^{\frac{S_{\tau,\chi}(z)-S_{\tau,\chi}(w)}{r}+O(1)}\frac{dzdw}{z-w}
\\+\frac{1}{2\pi \mathfrak{i}}\int_{\bar z_{cr}}^{z_{cr}}
\frac{\Phi_{\bar{q}}(z,t_1)}{\Phi_{\bar{q}}(z,t_2)}
z^{h_2-h_1-\frac 12 (t_1-t_2)-1}dz.
\end{multline*}
Since $\Re(S_{\tau,\chi}(z))\leq \Re(S_{\tau,\chi}(z_{cr}))\leq \Re(S_{\tau,\chi}(w))$ along the contours $C'_z$ and $C'_w$, the first integral converges to $0$ as $r\rightarrow 0$. Thus
\begin{equation*}
\lim_{r\rightarrow 0}K_{\bar{q}}((t_1,h_1),(t_2,h_2))
=
\lim_{r\rightarrow 0}\frac{1}{2\pi \mathfrak{i}}\int_{\bar z_{cr}}^{z_{cr}}
\frac{\Phi_{\bar{q}}(z,t_1)}{\Phi_{\bar{q}}(z,t_2)}
z^{h_2-h_1-\frac 12 (t_1-t_2)-1}dz.
\end{equation*}
Since $t_1,t_2\approx \frac{\tau}{r}>0$, from \eqref{eq:Phis} we have 
\begin{align}
\label{eq:PhiResidue}
\frac{\Phi_{\bar{q}}(z,t_1)}{\Phi_{\bar{q}}(z,t_2)}
&=\frac{\Phi^-_{\bar{q}}(z,t_1)}{\Phi^-_{\bar{q}}(z,t_2)}
\frac{\Phi^+_{\bar{q}}(z,t_2)}{\Phi^+_{\bar{q}}(z,t_1)}
=\frac{\Phi^+_{\bar{q}}(z,t_2)}{\Phi^+_{\bar{q}}(z,t_1)}
=\prod_{\substack{\min(t_1,t_2)<m<\max(t_1,t_2)\\ m\in \mathbb{Z}+\frac 12}}(1-zx_m^+)^{
\Delta t
},\\\nonumber
&=\prod_{\substack{\min(t_1,t_2)<m<\max(t_1,t_2)\\ m\in \mathbb{Z}+\frac 12}}\left(1-ze^{r(-m-\frac 12)}\beta^{-1}_{(d-m-\frac 12) \bmod k}\right)^{
\Delta t
}.
\end{align}

Let $N_{t,i}$ be the number of $t<m<d_r$ such that $\beta_{d-m-\frac 12}=\tilde\beta_i$ and let $\Delta N_i:=N_{t_1,i}-N_{t_2,i}$, which, up to a sign, counts the number of half-integers $m$ between $t_1$ and $t_2$ such that $\beta_{d-m-\frac 12}=\tilde\beta_i$.

Taking the limit $r\rightarrow 0$ we obtain the following theorem:
\begin{theorem}
\label{thm:bulkNearBoundary}
Let $r>0$. Consider plane partitions confined to an $\infty\times d_r$ box,  with periodic weights given by \eqref{eq:perwts} with $q=e^{-r}$. Let $\tilde\beta_i$'s be defined as in Definition \ref{def:betas}. Assume $d_r\bmod k$ is independent of $r$ and $\lim_{r\rightarrow 0}rd_r=V$. 
In the limit when $r$ goes to $0$ along a vertical line $\tau=V-\varepsilon$ near the edge of the system the liquid region consists of $l$ intervals, one for each $\tilde\beta_i$. The correlation functions of the system near a point $(\tau,\chi)$ in the bulk where $\chi$ is in the interval $(\chi_{j,+},\chi_{j,-})$ (see \eqref{eq:chiCr}) are given by determinants of the kernel
\begin{equation*}
\lim_{r\rightarrow 0}K_{\bar{q}}((t_1,h_1),(t_2,h_2))
=
\frac{1}{2\pi \mathfrak{i}}\int_{C}
\prod_{i=1}^l\left(1-ze^{-\tau}\tilde\beta^{-1}_{i}\right)^{\Delta N_i}
z^{-\Delta h-\frac 12 \Delta t-1}dz,
\end{equation*}
where $\Delta t=t_1-t_2$, $\Delta h=h_1-h_2$, $rt_i\rightarrow\tau$, $rh_i\rightarrow\chi$ and the integration contour $C$ connects the two non-real critical points of $S_{\tau,\chi}(z)$, passing through the real line in the interval $(\tilde\beta_l\gamma,e^\tau\tilde\beta_1)$ if $\Delta(t)\geq 0$ and through $(-\infty, 0)$ otherwise.
\end{theorem}

\begin{remark}
Notice, that $\sum_{i=1}^l \Delta N_i=\Delta t$. In particular, in the homogeneous case $\alpha_1=\dots=\alpha_k=1$, we have $l=1$ and $\tilde\beta_1=1$, so we recover the incomplete beta kernel. Otherwise, the process we have is translation invariant under translations by $k\mathbb{Z}\times\mathbb{Z}$ only and is a special case of processes studied previously in \cite{Bor} and \cite{BorShl2009}. 
\end{remark}

\begin{remark}
When $\chi$ is such that $S_{V-\varepsilon,\chi}$ has only real critical points, then during the deformation of the contours they do not cross each other and we do not pick up any residues from the $\frac 1{z-w}$ term. It follows that the correlation functions for horizontal lozenges decay exponentially fast as $r\rightarrow\infty$, so in the scaling limit almost surely there will be no horizontal lozenges near such a point $(V-\varepsilon,\chi)$. 
\end{remark}

\begin{remark}
The argument in Section \ref{sec:numCritPts} showing that $S_{V-\varepsilon,\chi}$ has only real critical points when $\chi$ is large also applies for $S_{0,\chi}$ without modification if the sets $\{\tilde{\beta}_i\}_{i=1}^l$ and $\{\tilde{\beta}_i\gamma\}_{i=1}^l$ are disjoint. Thus, in that case we also obtain that when $\tau=0$ and $\chi$ is large, almost surely the there will be no horizontal lozenges near the macroscopic point $(0,\chi)$. It follows that the height of the scaled plane partition is actually bounded, unlike the cases of homogeneous  weights studied in \cite{OR1},\cite{OR2},\cite{BMRT},\cite{M}.
\end{remark}

\subsubsection{Deformation of contours and the limiting correlation kernel near turning points} 
We now turn our attention to the point processes at the turning points. Let $rd_r=V$, $t_i=d_r-\hat t_i$, $h_i=\lfloor\frac{\chi_i}{r}\rfloor+\frac{\hat h_i}{r^{1/2}}$ with $\hat t_i,\hat h_i$ fixed, independent of $r$, and $\chi_j$ the height of the $j$'th turning point, which, from \eqref{eq:chiCr}, is
\begin{equation}
\label{eq:chiTurning}
\chi_{j}=-\frac V2-\frac1k\sum_{i=1}^lm_i\ln\left(\frac{e^V\tilde\beta_j-\tilde\beta_i\gamma}{e^V\tilde\beta_j}\right).
\end{equation}

From the earlier computations of critical points we see that the critical point corresponding to the turning point at height $\chi_i$ is $z_j=\tilde\beta_je^V$.

\begin{figure}[ht]
\caption{\label{fig:contourTurning}  The shaded region is the region $\Re(S_{V,\chi_j}(z))>\Re(S_{V,\chi_j}(z_{j}))$. The solid curves are the original contours whereas the dotted ones are the new, deformed contours.}
\includegraphics[width=13cm]{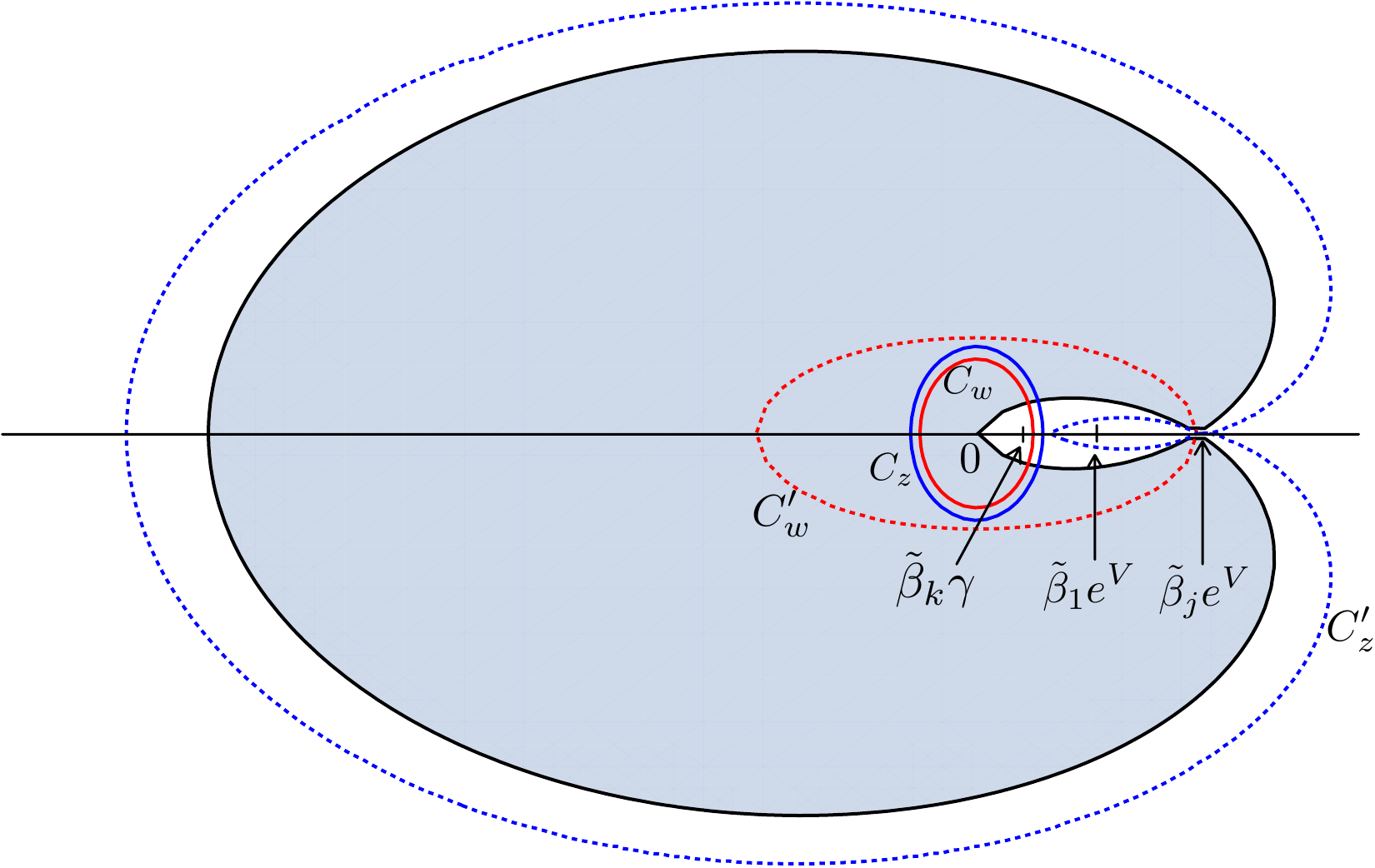}
\end{figure}

We have
\begin{equation*}
S_{\tau,\chi}(z)=\frac 1k\sum_{i=1}^k\int_{-\infty}^{0}\ln(1-z^{-1}\beta_i\gamma e^{M})dM
-\left(\chi+\frac V2 \right)\ln z,
\end{equation*}
\begin{equation*}
\left.\frac{d S_{V,\chi_j}(z)}{dz}\right|_{z=z_j}=0,
\end{equation*}
and
\begin{equation*}
\frac{d}{dz}\left(z\frac{d S_{V,\chi}(z)}{dz}\right)
=-\frac1k\sum_{i=1}^lm_i\left(\frac1{z-\tilde\beta_i\gamma}-\frac1{z}\right)
=-\frac1k\sum_{i=1}^l\frac{m_i\tilde\beta_i\gamma}{(z-\tilde\beta_i\gamma)z}.
\end{equation*}
It follows that 
\begin{equation*}
\frac{d}{dz}\left.\left(z\frac{d S_{V,\chi}(z)}{dz}\right)\right|_{z=z_j}<0,
\end{equation*}
which implies
$S''_{V,\chi}(z_j)<0$ for all $j$.

\begin{figure}[ht]
\caption{\label{fig:turningPointIntContour} The $\xi$ and $\omega$ contours.}
\includegraphics[width=2cm]{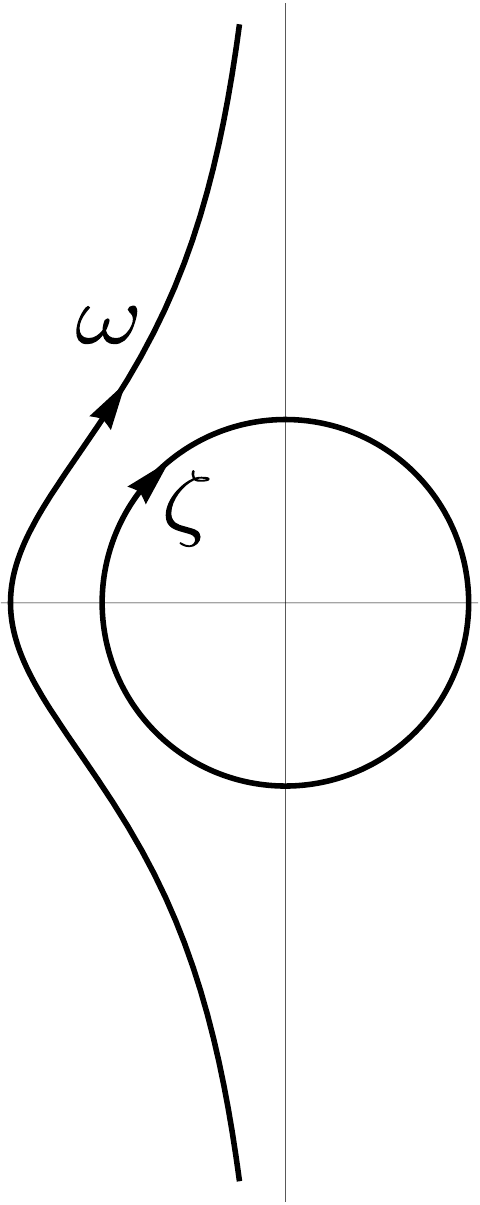}
\end{figure}

Now, the contours $\Re(S_{V,\chi_j}(z))=\Re(S_{V,\chi_j}(z_j))$ look as in Figure \ref{fig:contourTurning} and we deform the integration contours of the correlation kernel to pass through $z_j$ as shown in Figure \ref{fig:contourTurning}. When $\Delta t\geq 0$ the contours do not pass through each other so we do not pick up any resides from the $\frac 1{z-w}$ term. Along the new contours, away from the critical point we have $\Re(S_{V,\chi_j}(w))>\Re(S_{V,\chi_j}(z))$, so the contribution to the correlation kernel is exponentially small and the main contribution comes from the vicinity of the critical point. We make the change of variable $z=z_je^{r^{ 1/2}\xi}$, $w=z_je^{r^{1/2}\omega}$. In these coordinates the steepest descent contours are as in Figure \ref{fig:turningPointIntContour} (similarly to \cite{OR3}).

We have
\begin{equation*}
S_{V,\chi_j}(z)=S_{V,\chi_j}(z_j)+S'_{V,\chi_j}(z_j)(z-z_j)+\frac{S''_{V,\chi_j}(z_j)}2(z-z_j)^2+\cdots.
\end{equation*}
Since
$$z-z_j=-\sqrt{r}\xi z_j(1+O(\sqrt{r}),$$
we get
$$S_{V,\chi_j}(z)=S_{V,\chi_j}(z_j)+\frac{S''_{V,\chi_j}(z_j)}2r\xi^2z_j^2+O(r^{\frac 32}).$$

Using $z_j=\tilde\beta_je^V$, a similar computation to \eqref{eq:PhiResidue} gives
\begin{align*}
\Phi^+_{\bar{q}}(z,t_1)
&=\prod_{\substack{t_1<m<d_r\\ m\in \mathbb{Z}+\frac 12}}\left(1-\tilde\beta_je^Ve^{r^{1/2}\xi}e^{r(-m-\frac 12)}\beta^{-1}_{(d-m-\frac 12) \bmod k}\right).
\end{align*}
Since $t_1=d_r-\hat t_1$ with $\hat t_1$ fixed, this is a finite product. Moreover, $$e^Ve^{r(-m-\frac 12)}=e^{r(d_r-m-\frac 12)}=1+O(r),$$ from which it follows that
\begin{align*}
\Phi^+_{\bar{q}}(z,t_1)
&=\prod_{\substack{t_1<m<d_r\\ m\in \mathbb{Z}+\frac 12}}\left(1-\tilde\beta_je^{r^{1/2}\xi}\beta^{-1}_{(d-m-\frac 12) \bmod k}\right)(1+O(r))
\\&=(-r^{1/2}\xi)^{N_{t_1,j}}\prod_{\substack{i=1\\i\neq j}}^l(1-\tilde\beta_j\tilde\beta_i^{-1})^{N_{t_1,i}}(1+O(\sqrt r)).
\end{align*}
Thus, we obtain
\begin{multline*}
K_{\bar{q}}((t_1,h_1),(t_2,h_2))
=z_j^{-\Delta h-\frac\Delta t} 
(-r^{1/2})^{-\Delta N_j}
\prod_{\substack{i=1\\i\neq j}}^l(1-\tilde\beta_j\tilde\beta_i^{-1})^{-\Delta N_i}
\\
\times\frac{1}{(2\pi \mathfrak{i})^2}
\int_{z\in C'_z}\int_{w\in C'_w}
e^{\frac{S''_{V,\chi_j}(z_j)}2z_j^2(\xi^2-\omega^2)}
\frac{\omega^{N_{t_2,j}}}{\xi^{N_{t_1,j}}}
\frac{1}{\sqrt{r}(\xi-\omega)}\frac{e^{\hat h_2\omega}}{e^{\hat h_1\xi}}r(1+O(\sqrt r))\ d\xi\ d\omega.
\end{multline*}
When computing $\det K_{\bar q}((t_i,h_i),(t_j,h_j))_{i,j=1}^n$, the terms $$z_j^{-\Delta h-\frac\Delta t} 
(-r^{1/2})^{-\Delta N_j}
\prod_{\substack{i=1\\i\neq j}}^l(1-\tilde\beta_j\tilde\beta_i^{-1})^{-\Delta N_i}$$ cancel out. Thus, we obtain the following theorem.

\begin{theorem}
\label{thm:turningPoints}
Let $r>0$. Consider plane partitions confined to an $\infty\times d_r$ box,  with periodic weights given by \eqref{eq:perwts} with $q=e^{-r}$. Let $\tilde\beta_i$'s be defined as in Definition \ref{def:betas}. Assume $d_r\bmod k$ is independent of $r$ and $\lim_{r\rightarrow 0}rd_r=V$. 
In the limit when $r$ goes to $0$ near the vertical boundary $\tau=V$ the system develops $l$ turning points, one for each $\tilde\beta_i$. The $j$'th highest turning point is at position $(V,\chi_j)$, where $\chi_j$ is given by \eqref{eq:chiTurning}, and the correlation functions of the system near it are given by determinants of the kernel
\begin{equation}
\label{eq:kernelTurningPoint}
K^j_{\bar{q}}((\hat t_1,\hat h_1),(\hat t_2,\hat h_2))
=\frac{\sqrt{r}}{(2\pi \mathfrak{i})^2}
\int_{z\in C'_z}\int_{w\in C'_w}
e^{\frac{S''_{V,\chi_j}(z_j)}2z_j^2(\xi^2-\omega^2)}
\frac{\omega^{N_{\hat t_2,j}}}{\xi^{N_{\hat t_1,j}}}
\frac{e^{\hat h_2\omega}}{e^{\hat h_1\xi}}\frac{d\xi\ d\omega}{(\xi-\omega)},
\end{equation}
where $\hat t_i=d_r-t_i$ is the horizontal distance from the boundary, $\frac{\hat h_i}{r^{1/2}}$ is the vertical distance from $\chi_j$, $N_{t,i}$ is the number of $m\in\mathbb{Z}+\frac 12$ such that $d_r-t<m<d_r$ and $\beta_{d-m-\frac 12}=\tilde\beta_i$, and the integration contours are as in Figure \ref{fig:turningPointIntContour}.

\end{theorem}

\begin{remark}
When $\Delta t<0$, the $w$ contour crosses completely over the $z$ contour, so we pick up residues from the $\frac{1}{z-w}$ term along a simple contour going once around $0$. However, as explained in particular in \cite{OR2} or \cite[Section 5.0.3]{M2Per}, this does not alter the expression of the correlation kernel given in \eqref{eq:kernelTurningPoint}.
\end{remark}

\begin{remark}
While the limiting process is discrete in the horizontal coordinates $t$, it is continuous in the vertical coordinates $h$, so the right-hand side of \eqref{eq:kernelTurningPoint} converges to $0$ as $r\to 0$. In fact, if we condition on the vertical slices $\hat{t}_i$, then the point process in the limit is a continuous determinantal point process whose correlation measures have densities given by determinants $$\det \left(\frac{K^j_{\bar{q}}((\hat t_i,\hat h_i),(\hat t_j,\hat h_j))}{\sqrt{r}}\right)_{i,j}.$$
\end{remark}

\subsection{The GUE-corners process}
\label{sec:GUE-corners}
Fix $j$ with $1\leq j\leq l$ and choose integers $\bar{j}=(j_1,j_2,\dots)$ such that $N_{j_i,j}=i$ for all $i\in\mathbb{N}$. Note, that we must have $j_1<j_2<\cdots$. Given $\tilde{t_i}\in\mathbb{N}$ and $\tilde{h}_i\in\mathbb{R}$, let $\tilde{\rho_r}((\tilde t_1,\tilde h_1),\dots,(\tilde t_n,\tilde h_n))$ be the probability that there are horizontal lozenges at positions $(t_i,h_i)=(d_r-j_{\tilde t_i}, \lfloor\frac{\chi_i}{r}\rfloor+\frac{\tilde h_i}{r^{1/2}})$. It follows, that as $r\rightarrow 0$, we have 
$$\tilde{\rho_r}((\tilde t_1,\tilde h_1),\dots,(\tilde t_n,\tilde h_n))\rightarrow \det(K^{\bar j}_{\bar{q}}((\tilde t_i,\tilde h_i),(\tilde t_l,\tilde h_l)))_{i,l=1}^n,$$ with
\begin{equation}
\label{eq:kernelTurningPointGUE}
K^{\bar j}_{\bar{q}}((\tilde t_1,\tilde h_1),(\tilde t_2,\tilde h_2))
=\frac{\sqrt{r}}{(2\pi \mathfrak{i})^2}
\int_{z\in C'_z}\int_{w\in C'_w}
e^{\frac{S''_{V,\chi_j}(z_j)}2z_j^2(\xi^2-\omega^2)}
\frac{\omega^{\tilde t_2}}{\xi^{\tilde t_1}}
\frac{e^{\tilde h_2\omega}}{e^{\tilde h_1\xi}}\frac{d\xi\ d\omega}{(\xi-\omega)}.
\end{equation}
This is the correlation kernel for the GUE-corners process (see e.g. \cite{OR3}). In other words, if we restrict the point process at the $j$'th turning point to the slices $j_1,j_2,\dots$ then we get the GUE-corners process. In particular, since the $i$'th slice of the GUE-corners process has $i$ points, it follows that near the $j$'th turning point the $i$'th slice from the boundary has $N_{i,j}$ points.

We can describe the connection with the GUE-corners process differently. The number of horizontal lozenges on vertical slices near a turning point weakly increases as you move further away from the boundary. Pick any turning point and any collection of slices, such that the $i$'th slice has $i$ horizontal lozenges near the turning point. Then the point process near the turning point restricted to the selected slices will be the GUE-corners process. Note, that the GUE-corners process has $i$ particles on slice $i$, so to obtain the GUE-corners process it is obviously necessary to pick slices so that the $i$'th has $i$ particles. The key is that that's also sufficient, so as long as the number of particles is selected correctly, we obtain the GUE-corners process.

\subsection{The frozen regions separating the turning points}
\label{sec:FrozenRegions}
It follows from Theorems \ref{thm:bulkNearBoundary} and \ref{thm:turningPoints} that near the boundary at $\tau=V$ the system develops $l+1$ frozen regions separated by $l$ liquid regions. In the scaling limit the correlation functions for horizontal lozenges converge to zero in these frozen regions, so no horizontal lozenges appear in these regions. Let's call the lozenges of orientation 
\begin{tikzpicture}
	[scale=.3, thick]
	\lozl{(0,0)}{yellow!20!white}
\end{tikzpicture}
left lozenges and the lozenges of orientation 
\begin{tikzpicture}
	[scale=.3, thick]
	\lozr{(0,0)}{red!20!white}
\end{tikzpicture}
right lozenges. Since the frozen regions consist of only left and right lozenges, along any vertical line inside a frozen region only one type of lozenge can appear. Thus the frozen regions that appear can be characterized by specifying a horizontal cross-section, which can be done via a sequence of the letters $R$ and $L$ standing for right and left lozenges respectively. It is obvious that the top frozen region has only left lozenges, whereas the bottom frozen region has only right lozenges. 

First, let's consider the case when $l=k$, i.e. when $\beta_1,\dots,\beta_k$ are distinct, and let $\tilde\beta_j=\beta_{i_j}$ for $j=1,\dots,k$. From the discussion in Section \ref{sec:GUE-corners} it follows that the top turning point the first $i_1-1$ slices will have no horizontal lozenges, the next $k$ slices one horizontal lozenge each, the next $k$ slices two each, etc. It follows that the second frozen region from the top will have profile $ (L^{k-1}R)L^{i_1-1}$ where $*^m$ stands for $*$ repeated $m$ times, and $(*)$ stands for $*$ repeated indefinitely.

At the second turning point the first $i_2-1$ slices will have no horizontal lozenges, the next $k$ slices one horizontal lozenge each, the next $k$ slices two each, etc. It follows that the third frozen region from the top will have profile $ (L^{k-|i_1-i_2|-1}RL^{|i_1-i_2|-1}R)L^{\min(i_1,i_2)-1}$. It can inductively be shown, that after ``passing through'' a turning point, exactly one of the $L$'s in the repeating pattern of the frozen region changes to an $R$. 

Thus, all frozen regions near the right boundary are vertical facets. If we vertically project them to the floor and introduce polar coordinates so that the projection of the top facet $(L^k)$ has angle $0$ and the projection of the bottom facet $(R^k)$ has angle $\pi/2$, then the $i$'th facet from the top will project to angle $\frac{(i-1)\pi}{2k}$, and ignoring the first few slices, will be formed by a sequence of $i-1$ $R$'s and $k-i+1$ $L$'s repeated periodically. The particular pattern of $L$'s and $R$'s will depend on the order of $\beta_1,\dots,\beta_k$ and can be arbitrary.

\begin{figure}[ht]
\caption{\label{fig:4perFrozens} The frozen regions corresponding to the profiles \eqref{eq:4perprofiles}.}
\includegraphics[width=12cm]{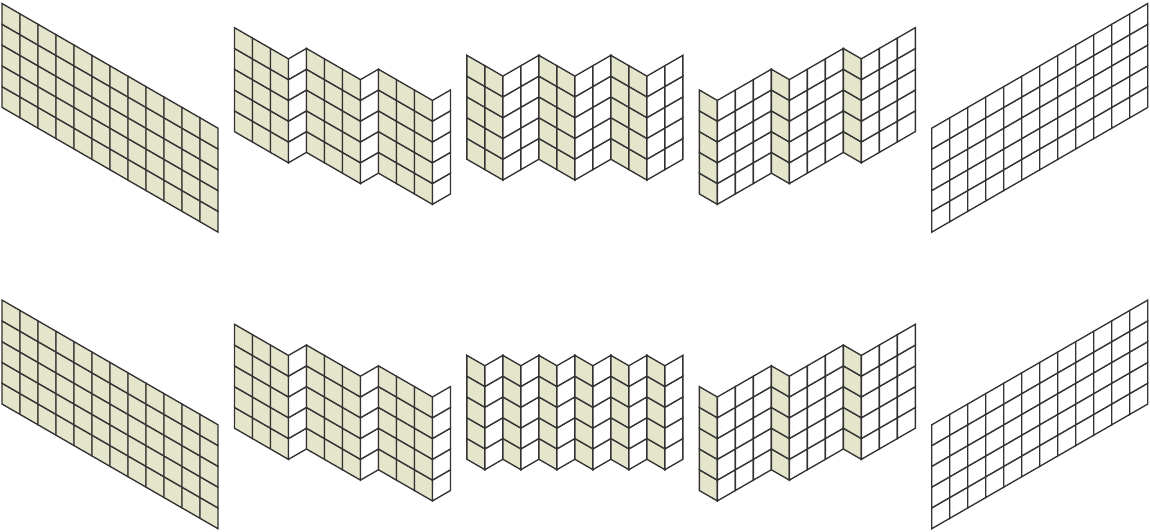}
\end{figure}

For example, for $k=4$, ignoring the first few slices, the frozen regions from top to bottom will have profiles
\begin{align}
\nonumber(L^4),(L^3R),&(L^2R^2),(LR^3),(R^4),
\\\label{eq:4perprofiles}&\text{or}
\\\nonumber(L^4),(L^3R),&(LRLR),(LR^3),(R^4),
\end{align}
depending on the order of the $\beta$'s (see Figure \ref{fig:4perFrozens}).

The $l<k$ case is similar, with the only difference being that some facets are missing. For example, if $k=4$ and $\beta_1=\beta_2>\beta_3>\beta_4$, then there will only be $3$ turning points and the $(L^3R)$ facet will be missing as the top two turning points will have merged.

\begin{figure}[ht]
\caption{\label{fig:simulations} An exact sample(generated using the algorithm in \cite{BorSchurDynam}) of plane partitions with $3$-periodic weights, modified at the slice through the corner. Near the right edge of the figure the $3$ turning points and the frozen regions separating them are visible. }
\includegraphics[width=12cm]{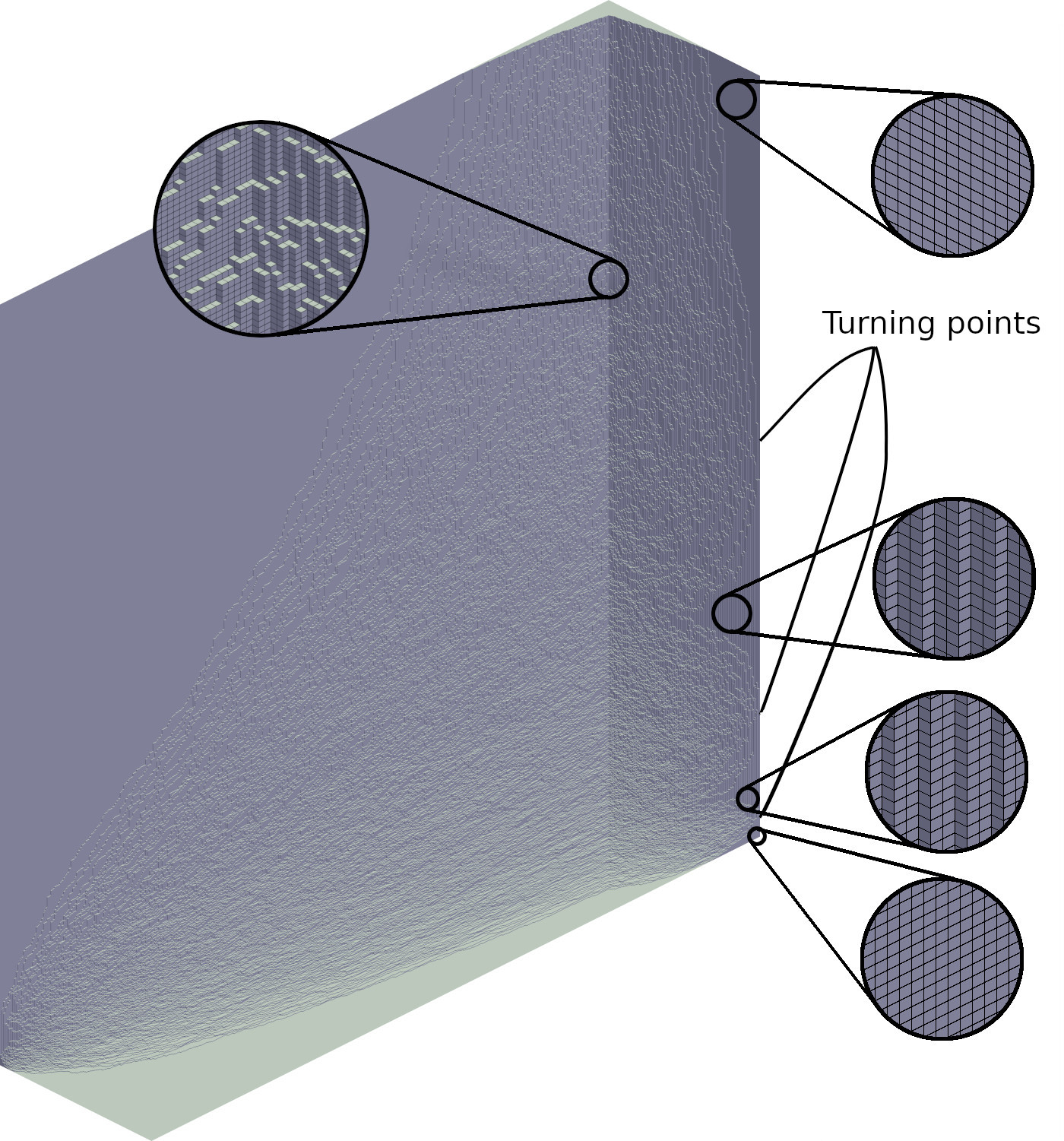}
\end{figure}

An exactly sampled plane partition with 3-periodic weights is shown in Figure \ref{fig:simulations}.

\section{Point processes at the first-order phase transition}
It can be seen from Figure \ref{fig:simulations} that the system develops a first order phase-transition at the slice $\tau=0$: the limit surface is not differentiable along $\tau=0$. The goal of this section is to study this phase transition, and obtain the point processes along it, which are not translation invariant in the horizontal direction. 

\begin{figure}[ht]
\caption{\label{fig:yms} The sequence $y_m$ for $m\in\mathbb{Z}+\frac12$.}
\includegraphics[width=14cm]{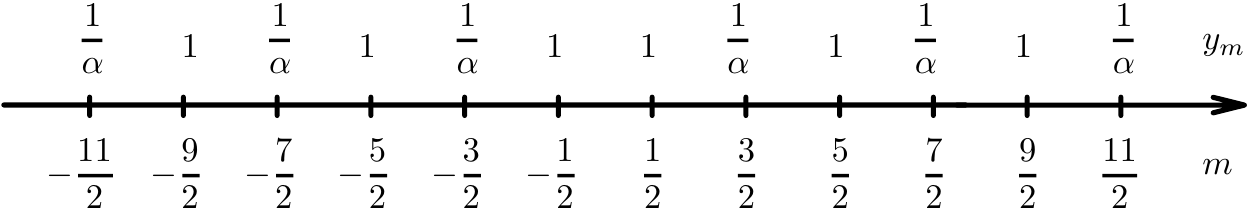}
\end{figure}

The nature of the fluctuations at $\tau=0$ should not depend on the precise boundary. To simplify the computations we set $V=\infty$ and restrict the derivations to the $k=2$ case. In this case the computations can be made explicitly. We have 
\begin{equation}
\label{eq:wtsformiddle}
\alpha_0=\alpha>1,\alpha_1=\frac 1\alpha,\gamma=\frac 1\alpha,\tilde\beta_1=1,\tilde\beta_2=\alpha.
\end{equation}
Since $d=\infty$, the form of $x_m^\pm$ from \eqref{eq:xs-prodqs} is not suitable anymore. Instead, we rewrite \eqref{eq:xs-qs} in the form
\begin{equation}
\label{eq:xs-middle}
x_m^+=q^my_m, x_{-m}^-=q^my_{-m}=q^my_m,\forall m\in\frac 12+\mathbb{Z}_{\geq 0},
\end{equation}
where $y_m$ is as in Figure \ref{fig:yms}. 

From \eqref{eq:Sp} and \eqref{eq:Spp} for $\tau\geq 0$ we have
\begin{equation}
\label{eq:Spt0}
2z\frac{dS_{\tau,\chi}}{dz}(z)=-\ln\left(\frac{z-1}{z}\right)-\ln\left(\frac{z-\frac 1\alpha}{z}\right)+\ln\left(\frac{-1}{z- e^\tau}\right)+\ln\left(\frac{-\alpha}{z-\alpha e^\tau}\right)-2\chi+\tau
\end{equation}
and
\begin{equation}
\label{eq:Sppt0}
2\frac{d}{dz}\left(z\frac{d S_{\tau,\chi}(z)}{dz}\right)
=-\frac{1}{z-1}
+\frac2z
-\frac1{z-\alpha^{-1}}
-\frac1{z-e^{\tau}\alpha}
-\frac1{z-e^{\tau}}.
\end{equation}

\subsubsection{The double real critical points}
As in Section \ref{sec:bulkandturningpoints} we first determine the points $(\tau,\chi)$ for which the function $S_{\tau,\chi}$ has double real critical points. Given $z\in\mathbb{R}$, we ask for which $(\tau,\chi)$ it is a double real critical point of $S_{\tau,\chi}$. If $z$ is a double real critical point of $S_{\tau,\chi}$, then $\tau$ can be obtained from 
\begin{equation}
\label{eq:Sppzero}
\frac{d}{dz}\left(z\frac{d S_{\tau,\chi}(z)}{dz}\right)=0
\end{equation}
and $\chi$ then can be easily computed from $z\frac{dS_{\tau,\chi}(z)}{dz}=0$. 

Suppose $\tau\geq 0$. Equation \eqref{eq:Sppzero} is quadratic in $e^\tau$. Solving it we obtain
\begin{equation}
\label{eq:DoubleCritTauPos}
e^\tau=z^2\text{, or }
e^\tau=\frac{z+\alpha z-2\alpha z^2}{2\alpha-\alpha z-\alpha^2 z}.
\end{equation}
%
%


If
\begin{equation}
\label{eq:etz}
e^\tau=\frac{z+\alpha z-2\alpha z^2}{2\alpha-\alpha z-\alpha^2 z},
\end{equation}
then $z>0$, since otherwise $e^\tau\leq 0$. It follows from \eqref{eq:Spt0} that $\chi$ is real if and only if
\begin{equation*}
1-z^{-1}>0,\ 1-\alpha^{-1}z^{-1}>0,\ 1-e^{-\tau}\alpha^{-1}z>0,\ \text{and }1-e^{-\tau}z>0,
\end{equation*}
or, equivalently, since $z>0$, that
\begin{equation*}
1<z<e^\tau=\frac{z+\alpha z-2\alpha z^2}{2\alpha-\alpha z-\alpha^2 z},
\end{equation*}
which, using \eqref{eq:etz}, becomes
\begin{equation*}
1<z\text{ and }1<\frac{1+\alpha-2\alpha z}{\alpha(2-z-\alpha z)}.
\end{equation*}
Since $\alpha>1$, we have $1<\alpha z$, which implies the denominator in the formula above is negative. Rewriting the second inequality we obtain
$(\alpha-1)(1-\alpha z)>0,$ which is a contradiction. Hence, if $(\tau,\chi)$ is a double real critical point and $\tau\geq 0$, we must have $e^\tau=z^2$. Similarly it can be established that when $\tau<0$, we still must have $e^\tau=z^2$. Thus, at the real critical points we have $z=\pm e^{\frac \tau 2}$.
\begin{figure}[ht]
\caption{\label{fig:chitau} A plot of the curve $\chi_\pm(\tau)$ when $a=9$.}
\includegraphics[width=8cm]{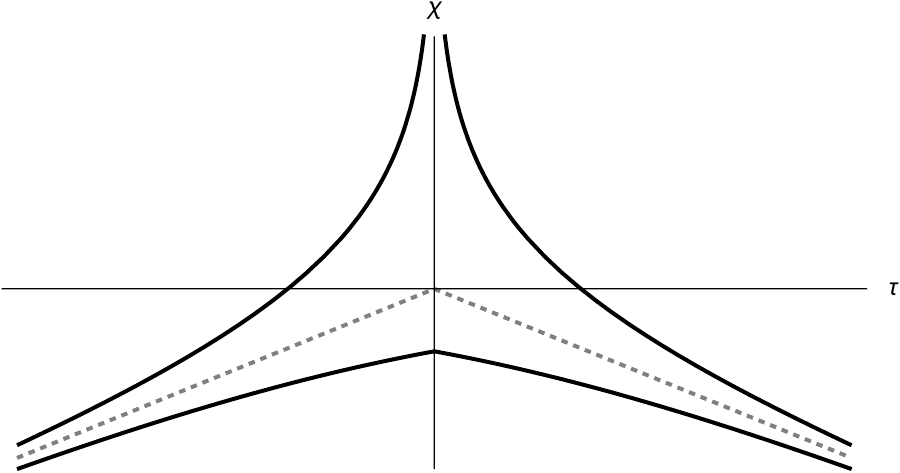}
\end{figure}
Solving $\frac{d S_{\tau,\chi}}{dz}(\pm e^{\frac\tau2})=0$ for $\chi$ we obtain that independently of the sign of $\tau$ we have
\begin{equation}
\label{eq:frozenBoundary-chi(tau)}
\chi_\pm(\tau)
=-\ln(1\mp e^{-\frac{|\tau|}2})
-\ln(1\mp \alpha^{-1}e^{-\frac{|\tau|}2})
-\frac12|\tau|.
\end{equation}
Figure \ref{fig:chitau} gives a plot of the curve $\chi_\pm(\tau)$. Note, that the curve $\chi_\pm(\tau)$ near $\tau=0$ is not differentiable. We have
\begin{equation*}
\frac{d\chi}{d\tau}
=\frac{\mp\frac{\sign(\tau)}2e^{-\frac{|\tau|}2}}{1\mp e^{-\frac{|\tau|}2}}
+\frac{\mp\frac{\sign(\tau)}2\alpha^{-1}e^{-\frac{|\tau|}2}}{1\mp \alpha^{-1}e^{-\frac{|\tau|}2}}
-\frac{\sign(\tau)}2.
\end{equation*}

If $z=e^{\frac\tau2}$, we have
\begin{equation*}
\lim_{\tau\rightarrow0\pm}\frac{d\chi}{d\tau}=\mp\infty.
\end{equation*}
If $z=-e^{\frac\tau2}$, we have
\begin{equation*}
\lim_{\tau\rightarrow0\pm}\frac{d\chi}{d\tau}=\pm\frac{1-\alpha}{4(1+\alpha)},
\end{equation*}
which implies that the frozen boundary is not differentiable at $\tau=0$.

\subsubsection{The number of critical points}
To carry out the steepest descent analysis, we need to know the number of complex non-real critical points of $S_{\tau,\chi}$. To obtain the number we exponentiate the equation
$2z\frac{d S_{\tau,\chi}(z)}{dz}=0$
to obtain
\begin{multline}
\label{eq:exponentiatedCritEqn}
\left(z-e^{\min\{\tau,0\}}\right)\left(z-e^{\min\{\tau,0\}}\alpha^{-1}\right)
\left(1-e^{-\max\{\tau,0\}}\alpha^{-1}z\right)
\left(1-e^{-\max\{\tau,0\}}z\right)
\\-e^{-2\left(\chi-\frac12V(\tau)\right)}z^2=0.
\end{multline}
The left-hand side is a degree $4$ polynomial with positive leading coefficient. Since its value is negative at $e^{\min\{\tau,0\}}$, the equation must have at least two distinct real solutions. Hence, the number of non-real complex solutions is at most $2$.

If we divide the $(\tau,\chi)$ plane into regions according to the number of critical points, $S_{\tau,\chi}$ will have double real critical points at the boundaries of those regions. Thus, the boundary is given by the curve $\chi_\pm(\tau)$ from \eqref{eq:frozenBoundary-chi(tau)}. It is easy to see that for any $\tau\neq 0$, if $\chi$ is very large, then the equation \eqref{eq:exponentiatedCritEqn} has four distinct real solutions, whence $S_{\tau,\chi}$ has no non-real critical points. It follows that in the region bounded by the curve $\chi_\pm(\tau)$ the function $S_{\tau,\chi}$ has a pair of complex conjugate critical points, whereas in the region outside the curve all critical points are real.

Carrying out the saddle point analysis is straightforward. When $\tau\neq 0$ we obtain that the limiting correlation functions at any point in the interior of the curve $\chi_\pm(\tau)$ are of the same type as in Theorem \ref{thm:bulkNearBoundary} and outside the curve we have frozen regions.

Let $\tau=0$ and $\chi>\chi_-(\tau)$. Then $S_{\tau,\chi}$ has exactly one pair of complex conjugate critical points. The deformation of contours and the saddle point analysis works exactly as in Theorem \ref{thm:bulkNearBoundary}. The only difference is in the residue terms picked up when the $z$ and $w$ contours cross each other. We have $t_1,t_2$  fixed, $\Delta t=t_1-t_2$ and $rh_i\rightarrow\chi$ with $\Delta h=h_1-h_2$ fixed as before. We obtain
\begin{equation*}
\lim_{r\rightarrow 0}K_{\bar{q}}((t_1,h_1),(t_2,h_2))
=
\lim_{r\rightarrow 0}
\frac{1}{2\pi \mathfrak{i}}\int_{C}
\frac{\Phi_{\bar{q}}(z,t_1)}{\Phi_{\bar{q}}(z,t_2)}
z^{-\Delta h-\frac 12 |t_1|+\frac 12|t_2|-1}dz,
\end{equation*}
where the contour $C$ is as before. From \eqref{eq:xs-middle} it follows that if $|m|<\max(|t_1|,|t_2|)$, then $x_m^\pm\rightarrow y_m$ as $r\rightarrow 0$. Given integers $t_1,t_2$, let $A(t_1,t_2)$ be the number of half integers $m$ between $t_1$ and $t_2$ such that $y_m=1$, taken with a positive sign if $t_1\geq t_2$ and with a negative sign otherwise. Define $g(t_1,t_2)$ by
\begin{equation*}
g(t_1,t_2)=
\begin{cases}
0,& t_1,t_2>0,\\
\Delta t,& t_1,t_2<0,\\
t_2,& t_2<0<t_1,\\
t_1,& t_1<0<t_2.
\end{cases}
\end{equation*}
We obtain
\begin{equation*}
\frac{\Phi_{\bar{q}}(z,t_1)}{\Phi_{\bar{q}}(z,t_2)}z^{-\Delta h-\frac 12 |t_1|+\frac 12|t_2|-1}
\rightarrow(-1)^{g(t_1,t_2)}(1-z)^{A(t_1,t_2)}(1-z/\alpha)^{\Delta t-A(t_1,t_2)}z^{-\Delta h-\frac {\Delta t}2-1}.
\end{equation*}
When taking determinants, the terms $(-1)^{g(t_1,t_2)}$ cancel out and we obtain the following result.
\begin{theorem}
\label{thm:bulkmiddle}
Let $r>0$. Consider plane partitions confined to an $\infty\times \infty$ box,  with periodic weights given by \eqref{eq:perwts},\eqref{eq:wtsformiddle} with $q=e^{-r}$. In the limit when $r$ goes to $0$ along the middle slice $\tau=0$ when $\chi>\chi_-(0)$, the correlation functions of the system (see \eqref{eq:chiCr}) are given by determinants of the kernel
\begin{equation*}
\lim_{r\rightarrow 0}K_{\bar{q}}((t_1,h_1),(t_2,h_2))
=
\frac{1}{2\pi \mathfrak{i}}\int_{C}
(1-z)^{A(t_1,t_2)}(1-z/\alpha)^{\Delta t-A(t_1,t_2)}z^{-\Delta h-\frac 12 \Delta t-1}dz,
\end{equation*}
where $\Delta t=t_1-t_2$, $\Delta h=h_1-h_2$, $t_1,t_2$ are fixed, $rh_i\rightarrow\chi$ with $\Delta h$ fixed, and the integration contour $C$ connects the two non-real critical points of $S_{0,\chi}(z)$, passing through the real line in the interval $(1,\alpha)$ if $\Delta(t)\geq 0$ and through $(-\infty, 0)$ otherwise.
\end{theorem}

\begin{remark}
Note, that this is a point process which is invariant under translations in the vertical direction but not in the horizontal direction. This is the case because the sequence $y_m$ is not translation invariant (see Figure \ref{fig:yms}). When $\alpha=1$ the process becomes the incomplete beta process from \cite{OR1}.
\end{remark}

\bibliography{mybib}
\bibliographystyle{alpha}

\end{document}